\documentclass[11pt]{amsart}
\usepackage{fullpage}
\usepackage{paralist} 
\usepackage{ifthen}
\usepackage{tropical}
\usepackage{eqnarray}
\usepackage{phdalgo}
\usepackage{array}
\usepackage{multirow}
\usepackage{url}
\usepackage{frenchineq}
\usepackage{tikz}

\usepackage{graphicx}
\usepackage{tikzexternal}
\tikzsetexternalprefix{figures/}
\usepackage{hyperref}

\renewcommand{\defi}{:=}
\newcommand{\etal}{\textit{et al.}}

\renewcommand{\mpray}[1]{\maxplus #1}

\newcommand{\BB}{\mathbb{B}}
\renewcommand{\PP}{\mathscr{P}}

\renewcommand{\CC}{\mathscr{C}}

\renewcommand{\HH}{\mathscr{H}}
\renewcommand{\DD}{\mathscr{D}}

\renewcommand{\tangent}{\mathscr{T}}
\newcommand{\GG}{\mathcal{G}}
\DeclareMathOperator{\closure}{\mathit{cl}}

\newcommand{\tper}{\mathop{\mathrm{tper}}}

\newcolumntype{M}[1]{>{\centering}m{#1}}
\newcolumntype{P}[1]{>{\centering}p{#1}}
\newcolumntype{B}[1]{>{\centering}b{#1}}

\newtheorem{definition}{Definition}
\newtheorem{proposition}{Proposition}
\newtheorem{propositiondefinition}[proposition]{Proposition-Definition}

\newtheorem{theorem}[proposition]{Theorem}
\newtheorem{corollary}[proposition]{Corollary}
\theoremstyle{remark}
\newtheorem{remark}{Remark}
\newtheorem{example}[remark]{Example}

\title{Computing the vertices of tropical polyhedra using directed hypergraphs}
\author{Xavier {A}llamigeon}
\address[X.~Allamigeon and S.~Gaubert]{INRIA and CMAP, \'Ecole Polytechnique, 91128 Palaiseau Cedex France}
\email[X.~Allamigeon]{xavier.allamigeon@inria.fr}
\author{{S}t{\'e}phane {G}aubert}
\email[S.~Gaubert]{stephane.gaubert@inria.fr}
\author{\'Eric Goubault}
\address[E.~Goubault]{CEA Saclay Nano-INNOV, Institut Carnot CEA LIST,  DILS/MeASI,
 Point Courrier n$^\circ$ 174, 91191 Gif sur Yvette Cedex, France}
\email{eric.goubault@cea.fr}
\date{\today}
\begin{document}
\begin{abstract}
We establish a characterization of the vertices of a tropical polyhedron defined  as the intersection of finitely many half-spaces. We show that a point is a vertex if, and only if, a directed hypergraph, constructed from the subdifferentials of the active constraints at this point, admits a unique strongly connected component that is maximal with respect to the reachability relation (all the other strongly connected components have access to it). This property can be checked in almost linear-time. This allows us to develop a tropical analogue of the classical double description method, which computes a minimal internal representation (in terms of vertices) of a polyhedron defined externally (by half-spaces or hyperplanes). We provide theoretical worst case complexity bounds and report extensive experimental tests performed using the library {\tt TPLib}, showing that this method outperforms the other existing approaches. 
\end{abstract}
\maketitle
\renewcommand{\thefootnote}{}\footnotetext{This work was partly performed when the first author was with EADS Innovation Works, SE/IS -- Suresnes, France and CEA, LIST MeASI -- Gif-sur-Yvette, France.}\footnotetext{The authors were partially supported by the Arpege programme of the French National Agency of Research (ANR), project ``ASOPT'', number ANR-08-SEGI-005. The last two authors were also partially supported by the Digiteo project DIM08 ``PASO'' number 3389.}
\renewcommand{\thefootnote}{\arabic{footnote}}

\section{Introduction}\label{sec:introduction}

Tropical polyhedra are the analogues of convex polyhedra in tropical algebra. The latter deals with structures like the max-plus semiring, which is the set $\real \cup \{-\infty\}$, equipped with the addition $(x,y)\mapsto \max(x,y)$ and the multiplication $(x,y)\mapsto x+y$.

The study of the tropical analogues of convex sets is an active research topic, which has been treated under various 
aspects. It arose in the work of Zimmermann~\cite{zimmermann77}, following a way opened by Vorobyev~\cite{vorobyev67}, motivated by optimization theory. Max-plus convex cones, thought of as the analogues of linear spaces, were studied by Cuninghame-Green~\cite{CG}. Their theory was independently developed by Litvinov, Maslov, and Shpiz (see in particular~\cite{litvinov00}, and also~\cite{maslov92}) 
with motivations from variations calculus and asymptotic
analysis, and by Cohen, Gaubert, and Quadrat~\cite{cgq00,cgq02} (see also~\cite{maxplus97}) who initiated a ``geometric approach'' of discrete event systems~\cite{ccggq99}, further developed by Katz~\cite{katz05,loiseau}. In~\cite{cgqs05,NiticaSinger07a}, Singer, Nitica, and some of the aforementioned authors, related this theory to abstract convexity~\cite{ACA}. The work of Briec and Horvath~\cite{BriecHorvath04} is also in the spirit of generalized convexity, some motivations from mathematical economy appeared in~\cite{BriecHorvath2009}.
Polyhedral max-plus convex sets also appeared in the work of Bezem, Nieuwenhuis, and Rodr\'{\i}guez-Carbonel~\cite{BezemNieuwenhuisCarbonnell10}, as sets defined by ``max-atoms'', with motivations from SMT (sat-modulo theory) solving.
Moreover, the field has been considerably developed after the work of Develin and Sturmfels~\cite{DS}, who related tropical and discrete geometry, showing in particular that tropical polyhedra can be thought of as regular polyhedral subdivisions of the products of two simplices. This was at the origin of a number of works,
by Joswig, Santos, Yu, Block, Ardila, and the same authors~\cite{joswig04,DSS,DevelinYu,blockyu06,JSY07,Joswig2008,ardila}. 

From the perspective of tropical geometry, tropical polyhedra may be thought
of as degenerate limits of classical polyhedra along a logarithmic deformation
(see~\cite{BriecHorvath04} for a proof of this fact), or as the image by the valuation of polyhedra over an ordered field of real Puiseux series. This explains a certain analogy between tropical and classical convexity. 
In particular, tropical analogues of several theorems in classical convexity have been established, including the ones of Hahn and Banach~\cite{zimmermann77,cgqs05,DS}, 
Minkowski~\cite{GK,BSS}, Minkowski-Weyl~\cite{GK06a,SGKatz11}, Radon~\cite{But,GM08}, Helly and Carath\'eodory~\cite{BriecHorvath04,gauser,GM08}, and also more advanced discrete convexity results~\cite{GM08}.

In contrast, algorithmic aspects of tropical polyhedra have not yet been thoroughly explored. In particular, a tropical polyhedron can be represented in two different ways, either externally, in terms of affine inequalities, or internally, as a set generated by finitely many points and rays, 
see~\cite{SGKatz11} and the references therein. The minimal internal representations of a tropical polyhedron are essentially unique, and consists of its extreme points (vertices) and representatives of extreme rays. Passing from an external description of a polyhedron to a (minimal) internal description, or inversely, is a fundamental computational issue, comparable to the well-known vertex/facet enumeration or convex hull problems in the classical case. 

In the present paper, we develop a combinatorial characterization of the extreme points and rays of tropical polyhedra defined externally. The characterization is equivalently expressed in terms of tropical polyhedral cones (as homogeneous representations of polyhedra). Polyhedral cones are sets consisting of vectors $\vect{x} = (\vect{x}_1,\ldots,\vect{x}_d)$ 
with entries in 
$\real \cup\{-\infty\}$ satisfying a system of linear inequalities in the tropical sense, \ie{} of the form:
\begin{equation}\label{e-fond}
\max_{i \in \oneto{d}} A_{ki}+\vect{x}_i \leq \max_{j \in \oneto{d}} B_{kj}+\vect{x}_j, \qquad \text{for } k \in \oneto{p},
\end{equation}
where for all integers $n$, $\oneto{n}$ refers to the set $\{1, \dots, n\}$, and $A, B$ are matrices of size $p \times d$ with entries in $\real \cup\{-\infty\}$.
If $\CC$ refers to the cone defined by the latter inequalities, a vector $\vect{v} \in \CC$ is said to be (tropically) extreme if it cannot be written as the point-wise supremum of two vectors of $\CC$ that are both different from it. We denote by $\argmax A_k \vect{v}$ (resp.\ $\argmax B_k \vect{v}$) the
set of indices $j \in \oneto{d}$ attaining the maximum at the left-hand side
(resp.\ right-hand side) of each inequality~\eqref{e-fond}. We associate with a vector
$\vect{v} \in \CC$ a directed hypergraph, referred to as the \emph{tangent directed hypergraph} at $\vect{v}$ in the cone $\CC$, consisting of the nodes $\{i \in \oneto{d} \mid \vect{v}_i \neq -\infty\}$, and 
one directed hyperarc $(\argmax B_k \vect{v} ,\argmax A_k \vect{v})$
for each index $k \in \oneto{p}$ such that both maxima in~\eqref{e-fond} coincide and take a finite value. This definition is illustrated in Section~\ref{sec:combinatorial_characterization}, in which more information on directed hypergraphs
can be found. The reachability relation induces a partial order on the strongly connected components of a directed hypergraph, meaning that a component is ``greater'' than another if the former can be reached from the latter. The main result of this paper is the following characterization: 
\begin{theorem}\label{th:main}\label{th-main}
Let $\CC$ be a tropical polyhedral cone. A vector $\vect{v} \in \CC$ is tropically extreme if, and only if, the set of the strongly connected components of the tangent directed hypergraph at $\vect{v}$ in $\CC$, partially ordered by the reachability relation, admits a greatest element.
\end{theorem}

This theorem shows interesting analogies and discrepancies with
the classical result stating that a point of a polyhedron defined by inequality constraints is a vertex if, and only if, the family of gradients of active constraints at this point is of full rank. In the tropical case, the expressions arising on both side of the constraints~\eqref{e-fond} are not differentiable, but they are convex, and so, they admit a subdifferential at each point at which they take a finite value. The subdifferential of the map $\vect{x} \mapsto \max_{i \in \oneto{d}} A_{ki}+\vect{x}_i$ at point $\vect{v}$ is easily seen to be the convex hull of the set of vectors of the canonical basis of $\real^d$ with indices in $\argmax A_k \vect{v}$. The same is true, mutatis mutandis, for the map appearing at the right-hand side of~\eqref{e-fond}. Hence, Theorem~\ref{th-main} appears to be an infinitesimal characterization, as the classical result. However, the classical rank condition does not have a tropical analogue: several rank notions have been considered in the tropical setting~\cite{DSS,AGG08b}, none of which explains the reachability condition appearing in Theorem~\ref{th-main}. 

Theorem~\ref{th-main} has both theoretical and algorithmic applications. In the companion paper~\cite{AllamigeonGaubertKatzJCTA2011}, it is used to show that 
the tropical analogues of the polar of the cyclic polytope have fewer
vertices than in the classical case (in other words, along the
deformation sending a classical polyhedron to a tropical polyhedron, some
classical extreme points degenerate in points which are no longer extreme
in the tropical sense). 

From the algorithmic point of view, a significant advantage of the criterion provided by Theorem~\ref{th-main} is that it can evaluated in almost linear time in the size of the tangent hypergraph (Theorem~\ref{th:complexity}). Thus, the corresponding computational complexity exclusively depends on the size of the external representation of the cone. This allows us to define an algorithm determining the extreme points and rays of a tropical polyhedron defined by inequalities (Section~\ref{sec:external_to_internal}). We call this algorithm the \emph{tropical double description method}, by analogy with the classical method which goes back to Motzkin~\etal{}~\cite{MRTT53} and was later refined by Fukuda and Prodon~\cite{FukudaProdon96}. Given a polyhedron defined by a system of $p$ inequalities, it consists in determining the set of the extreme generators of the polyhedron defined by the first $k$ inequalities, by induction on $k = 1, \dots, p$. It is based on a result (Theorem~\ref{th:ddm_basis}) allowing to build a set of generators of the intersection of a polyhedron with a half-space. This result can be extended to the intersection with tropical hyperplanes (Theorems~\ref{th:ddm_basis_hyperplane} and~\ref{th:ddm_basis_eq}), so that the tropical double description method can also handle polyhedra defined as mixed intersections of half-spaces and hyperplanes. Theorem~\ref{th-main} is the cornerstone of the double description method, since the latter algorithm critically relies on an efficient criterion to eliminate non-extreme generators (propagating such generators in the induction considerably increases the time complexity).

We include for the sake of comparison an alternative algorithm (Section~\ref{sec:arrangement}), based on determining the extreme generators of a polyhedron $\PP$ by computing the vertices of the arrangement formed by (tropical) hyperplanes associated with the half-spaces defining $\PP$, assuming that they are in general position. For some polyhedra, this algorithm has a better worst-case complexity than the double description method. However, its interest is rather theoretical, since this worst-case complexity is essentially tight, and it does not apply to arbitrary polyhedra. 

The inductive approach used in the tropical double description presented here is reminiscent of an algorithm of Butkovi\v{c} and Hegedus~\cite{butkovicH} computing a generating set of a tropical polyhedral cone described by linear (in)equalities. Gaubert gave a similar one and derived the equivalence between the internal and external representations~\cite[Ch.~III]{gaubert92a} (see also~\cite{maxplus97,SGKatz11}). Our approach is more general in the sense that it handles intersections with other kinds of constraints. Moreover, the efficient elimination of redundant candidates using directed hypergraphs brings an important breakthrough both in theory and in practice in comparison with the previous techniques. We refer the reader to Section~\ref{sec:comparison} for an exhaustive discussion.

In~\cite{Joswig2008}, Joswig defined a method which is able to compute the vertices of the polyhedral complex associated with a tropical polytope (in the sense of~\cite{DS}), from a set of generating points. Other approaches~\cite{LorenzoPuente2011,Truffet2010} rely on a similar technique applying on cones described by (in)equalities. While such algorithms are of interest from a combinatorial point of view, the size of the complex may be much larger than the number of vertices, leading to a suboptimal method to determine concise internal representations.

The dual problem of computing an external representation of a tropical polyhedron generated by a set of points and rays recently appeared
to be more tractable: in a paper of the first two authors with Katz~\cite{AllamigeonGaubertKatzLAA2011}, it is shown that such a representation can be determined in incremental quasi-polynomial time, hence with a total complexity quasi-polynomial in the size of the input and output. 
This result is based on the particular structure of polar cones of tropical polyhedra (relations between defining inequalities and weighted transversals in undirected
hypergraphs) and it cannot be transposed to the primal problem discussed in the present work. 

We also note that the tropical double description method 
allows one in particular to check whether 
the intersection of a family of half-spaces is empty. However, if one is only interested in checking the latter emptyness
property, different algorithms may be used. Indeed, the emptyness problem
is equivalent~\cite{AGGut09} to solving a mean payoff
game, a problem for which several combinatorial algorithms
have been developed, including pseudo-polynomial algorithms (no polynomial time algorithm is currently known).

Finally, we note that the main results of this paper have been announced in the proceedings article~\cite{AGG10}.

\section{Preliminaries on tropical polyhedra and polyhedral cones}\label{sec:preliminaries}

We denote by $\maxplus\defi  \real \cup \{-\infty\}$ the tropical (max-plus)
semiring. It is equipped with the addition $x \mpplus y \defi \max(x,y)$ and the multiplication $x \mptimes y \defi x + y$ (also denoted by concatenation $x y$). The neutral elements for these two laws are denoted by $\mpzero \defi -\infty$ and $\mpone \defi 0$. We shall use the notation $\lambda^{-1}$ for the tropical inverse of a scalar $\lambda\in \maxplus\setminus\{\mpzero\}$, which is nothing but the opposite of $\lambda$. 

The set $\maxplus^d$ refers to the $d$-th fold Cartesian
product of the tropical semiring. Its elements can be thought of as points of an affine space, or as vectors. They are denoted by bold symbols, for instance $\vect{x} = (\vect{x}_1,\dots,\vect{x}_d)$. 
The elements $\mpzerovect$ and $\mponevect$ refer to the vectors whose coordinates are all equal to $\mpzero$ and $\mpone$ respectively. In the sequel,
the tropical semiring $\maxplus$ will be equipped with the topology arising from the metric $(s,t)\mapsto |e^s-e^t|$. The set $\maxplus^d$ will be equipped with the product topology, and the associated closure operator will be denoted by $\closure(\cdot)$.

Tropical operations are naturally extended to vectors and matrices over $\maxplus$, defining $(A \mpplus B)_{ij} = A_{ij} \mpplus B_{ij}$ and $(A B)_{ij} = \mptimes_k A_{ik} B_{kj}$. 
The \emph{Minkowski sum} of two sets $S, S' \subset \maxplus^d$, denoted by $S \mpplus S'$, is defined as the set $\{ \vect{x} \mpplus \vect{x}' \mid (\vect{x}, \vect{x}') \in S \times S' \}$.

A set $\CC \subset \maxplus^d$ is said to be a \emph{tropical convex set} if it contains the \emph{tropical segments} between any two of its points $\vect{x}$ and $\vect{y}$. The latter is defined as the set of points of the form $\lambda \vect{x} \mpplus \mu \vect{y}$, for $\lambda, \mu \in \maxplus$ such that $\lambda \mpplus \mu = \mpone$. Note that this definition is analogous to the familiar one
(in standard convexity) which requires in addition
the scalars $\lambda$ and $\mu$ to be nonnegative: the latter
condition is automatically satisfied in the tropical setting,
since $\mpzero =-\infty\leq \lambda$ holds for any $\lambda \in \maxplus$. The \emph{tropical convex hull} $\mpco(S)$ of a subset $S \subset \maxplus^d$ is the set of the combinations $\lambda_1 \vect{x}_1 \mpplus \cdots \mpplus \lambda_p \vect{x}_p$, where $p \geq 1$, $\vect{x}_i \in S$ and $\lambda_i \in \maxplus$ for all $i\in[p]$, and $\lambda_1 \mpplus \cdots \mpplus \lambda_p = \mpone$.

Similarly, a set $\CC \subset \maxplus^d$ is said to be a \emph{tropical (convex) cone} if it contains all the combinations $\lambda \vect{x} \mpplus \mu \vect{y}$ ($\lambda, \mu \in \maxplus$) of any of two elements $\vect{x},\vect{y} \in \CC$. Given $S \subset \maxplus^d$, the \emph{tropical cone generated by $S$}, denoted by $\mpcone(S)$, is the set of the elements $\lambda_1 \vect{x}_1 \mpplus \cdots \mpplus \lambda_p \vect{x}_p$ where $p \geq 1$, $\vect{x}_i \in S$ and $\lambda_i \in \maxplus$ for all $i\in[p]$. A tropical convex set (resp.\ convex cone) 
is said to be
\emph{finitely generated} if it is of the form $\mpco(S)$ 
(resp.\ $\mpcone(S)$) for some
finite subset $S \subset \maxplus^d$. 
In the sequel, the terms \emph{convex set} or \emph{cone} are interpreted in the tropical sense.

Given a convex set $\CC \subset \maxplus^d$, a point $\vect{p} \in \CC$ is said to be an \emph{extreme point} (or \emph{vertex}) of $\CC$ if for all $\vect{x}, \vect{y} \in \CC$ and $\lambda,\mu \in \maxplus$ such that $\lambda \mpplus \mu = \mpone$, $\vect{p} = \lambda \vect{x} \mpplus \mu \vect{y}$ holds only if $\vect{p} = \vect{x}$ or $\vect{p} = \vect{y}$. Analogously, when $\CC$ is a convex cone, a non-null vector $\vect{v} \in \CC$ is said to be \emph{extreme} in $\CC$ if for all $\vect{x}, \vect{y} \in \CC$, $\vect{v} = \vect{x} \mpplus \vect{y}$ implies $\vect{v} = \vect{x}$ or $\vect{v} = \vect{y}$. In this case, the set $\mpray{\vect{v}} = \{ \lambda \vect{v} \mid \lambda \in \maxplus \}$ is said to be an \emph{extreme ray} of $\CC$, and the vector $\vect{v}$ is a \emph{representative} of this ray.

The subset $\BB \defi \{\mpzero, \mpone\}$ of the tropical semiring constitutes a sub-semiring of $\maxplus$. A subset $\DD$ of $\BB^d$ is said to be a \emph{boolean cone} if for all $\vect{x}, \vect{y} \in \DD$ and $\lambda, \mu \in \BB$, $\lambda \vect{x} \mpplus \mu \vect{y} \in \DD$ (in other words, $\DD$ is a sup-semilattice for the standard partial order on $\BB^d$). A vector of a boolean cone $\DD$ is said to be \emph{extreme} if it cannot be expressed as the pointwise supremum of two other vectors of $\DD$.

Tropical polyhedra and polyhedral cones are defined analogously to classical ones. A \emph{tropical affine half-space} is a set formed by the solutions $\vect{x} = (\vect{x}_i) \in \maxplus^d$ of a tropical affine inequality
\[
a_0 \mpplus  \mpsum_{i\in[d]} a_i \vect{x}_i \leq b_0 \mpplus \mpsum_{i\in[d]} b_i \vect{x}_i
\]
where $a_i, b_i \in \maxplus$ for all $i = 0,\dots,d$. It is said to be a \emph{tropical (linear) half-space} when the coefficients $a_0$ and $b_0$ are omitted. 
In this setting, a \emph{tropical polyhedron} (resp.\ a \emph{tropical polyhedral cone}) is the intersection of finitely many tropical affine (resp.\ linear) half-spaces. Equivalently, any tropical polyhedron can be seen as the set of the solutions of a system of inequality constraints $A \vect{x} \mpplus \vect{c} \leq B \vect{x} \mpplus \vect{d}$, where $A$ and $B$ are $p \times d$-matrices with entries in $\maxplus$, $\vect{c}$ and $\vect{d}$ are vectors of $\maxplus^p$, and $\leq$ denotes the standard partial ordering of vectors. Similarly, a tropical polyhedral cone is the set of the solutions of a two-sided system of the form $A \vect{x} \leq B \vect{x}$.

The description of tropical polyhedra and polyhedral cones as intersections of half-spaces 
is said to be \emph{external}. Moreover, tropical polyhedra and polyhedral cones admit an \emph{internal} representation,
by means of finitely many points and rays, as established by the following tropical analog of the Minkowski-Weyl theorem.
\begin{theorem}[{\cite[Th.~2]{SGKatz11}}]\label{th:minkowski_weyl}
The tropical polyhedra of $\maxplus^d$ are precisely the sets of the form $\mpco(P) \mpplus \mpcone(R)$ where $P$ and $R$ are finite subsets of $\maxplus^d$.

The tropical polyhedral cones of $\maxplus^d$ are precisely the sets of the form $\mpcone(V)$ where $V$ is a finite subset of $\maxplus^d$.
\end{theorem}
Thus, a tropical polyhedron $\PP$ is the sum of a bounded (finitely generated and convex) set and of a polyhedral cone. The latter coincides with the \emph{recession cone} $\mprec(\PP)$ of $\PP$, which is defined as the set $\{ \vect{v} \mid \vect{x} \mpplus \lambda \vect{v} \in \PP \text{ for all } \lambda \in \maxplus \}$, given an arbitrary point $\vect{x} \in \PP$, see~\cite{GK}.
 When $\PP \neq \emptyset$ is defined by a system of inequalities $A \vect{x} \mpplus \vect{c} \leq B \vect{x} \mpplus \vect{d}$, the recession cone can be shown to be the set of the solutions of the system $A \vect{x} \leq B \vect{x}$.

The couple $(P,R)$ is said to be a \emph{generating representation} of a tropical polyhedron $\PP$ when $\PP = \mpco(P) \mpplus \mpcone(R)$. Similarly, the set $V$ is a \emph{generating set} of a tropical polyhedral cone $\CC$ when $\CC = \mpcone(V)$. 
For algorithmic purposes, we look for representations which are
minimal in a suitable sense.
The following proposition, which combines several results of~\cite{GK},
shows that such representations do exist. (Actually, the results of~\cite{GK} apply more generally to \emph{closed} ---not necessarily polyhedral--- tropical convex sets.)
\begin{proposition}[{\cite[Th.~3.2, Th.~3.3, Coro.~3.4]{GK}}]\label{prop:minimal_gen}
A tropical polyhedron $\PP \subset \maxplus^d$ admits
a generating representation $(P,R)$ in which $P$ consists
of the extreme points of $\PP$, and $R$ contains precisely
one representative of each extreme ray of the recession cone
of $\PP$. Moreover, if $(P',R')$ is any generating representation
of $\PP$, then $P'\supset P$, and $R'$ contains at least one
scalar multiple of every element of $R$.

Similarly, a tropical polyhedral cone $\CC\subset \maxplus^d$ 
admits a generating set $V$ consisting of precisely one element
in each extreme ray of $\CC$. Moreover, if $V'$ is any generating
set of $\CC$, then $V'$ contains at least one scalar multiple
of every element of $V$.

\end{proposition}
The generating representations $(P,R)$ and the generating set
$V$ arising in this proposition will be referred to as \emph{minimal}.
The minimal generating representations of a tropical polyhedron (or of a tropical polyhedral cone) are almost identical, since they only differ by multiplicative factors on the representatives of extreme rays. We obtain \emph{canonical} minimal representations by requiring these vectors to be \emph{scaled} for the ``norm'' $\mpnorm{\cdot}$ over $\maxplus^d$ defined by $\mpnorm{\vect{x}} \defi \max_{i\in[d]} e^{\vect{x}_i}$, \ie{} to satisfy $\mpnorm{\vect{x}} = 1$. 

Tropical polyhedra of $\maxplus^d$ can be represented by polyhedral cones of $\maxplus^{d+1}$. In the classical setting, such a technique is known as \emph{homogenization} (see for instance Ziegler's monograph~\cite{ziegler98}). 
As shown in~\cite{cgq02,GK}, the same technique works in the tropical setting.

We restrict here our attention to the case of finitely generated convex sets. The notation $(M \ \vect{v})$ refers to the matrix obtained by appending the vector $\vect{v}$ after the last column of the matrix $M$.
\begin{definition}\label{def:homogenization}
Let $\PP  = \{ \vect{x} \in \maxplus^d \mid A \vect{x} \mpplus \vect{c} \leq B \vect{x} \mpplus \vect{d} \}$ be a non-empty tropical polyhedron ($A, B \in \maxplus^{p \times d}$, $\vect{c},\vect{d} \in \maxplus^d$). The \emph{homogenized cone} $\homo{\PP}$ is the polyhedral cone given by:
\[
\homo{\PP} \defi \bigl\{\vect{x} \in \maxplus^{d+1} \mid
(A \ \vect{c}) \, \vect{x} \leq 
(B \ \vect{d}) \, \vect{x} \bigr\}.
\]
\end{definition}

When $\vect{x} \in \maxplus^d$ and $\alpha \in \maxplus$, the element $(\vect{x},\alpha)$ refers to the vector of $\maxplus^{d+1}$ whose $d$ first coordinates coincide with $\vect{x}$, and the last coordinate is equal to $\alpha$. 
The generating representations of a tropical polyhedron and of its homogenized cone are connected by the following result, which is an immediate consequence
of the relations between a convex set and its homogenized cone which
are established in~\cite[\S~2]{GK}.
\begin{proposition}[Coro.\ of \cite{GK}]\label{prop:homo_gen}
Let $\PP \subset \maxplus^d$ be a non-empty tropical polyhedron. Then the following statements hold:
\begin{enumerate}[(i)]
\item if $(P,R)$ is a generating representation of $\PP$, then $(P \times \{\mpone\}) \cup (R \times \{\mpzero\})$ is a generating set of its homogenized cone $\homo{\PP}$.

\item conversely, if $V$ is a generating set of $\homo{\PP}$, then the couple $(P,R)$ defined by $P \defi \{ \alpha^{-1} \vect{p} \mid (\vect{p}, \alpha) \in V \text{ and } \alpha \neq \mpzero \}$ and $R \defi \{ \vect{r} \mid (\vect{r}, \mpzero) \in V \}$ forms a generating representation of $\PP$.

\item in the two previous statements, if any of the representations is minimal (and canonical), then the other is also minimal (and canonical).
\end{enumerate}
\end{proposition}
As a consequence, $\vect{p}$ is a vertex of $\PP$ if, and only if, the vector $(\vect{p}, \mpone)$ is an extreme vector of the homogenized cone $\homo{\PP}$. Similarly, the extreme vectors of the recession cone $\mprec(\PP)$ are precisely the elements $\vect{r} \in \maxplus^d$ such that $(\vect{r}, \mpzero)$ is extreme in $\homo{\PP}$. Thus, we will only state the main results of this work for tropical cones, leaving to the reader the derivation of the affine analogues using homogenization, along the lines of Proposition~\ref{prop:homo_gen}.

\begin{figure}[t]
\begin{center}
\begin{tikzpicture}
\begin{scope}[scale=0.8,>=triangle 45,convex/.style={draw=lightgray,fill=lightgray,fill opacity=0.7},point/.style={blue!50},line/.style={blue!50,ultra thick},convexborder/.style={ultra thick},pointlabel/.style={black},ray/.style={ultra thick,green!60!black}]
\draw[gray!30,very thin] (-3.5,-2.5) grid (3.5,4.5);
\draw[gray!50,->] (-3.5,0) -- (3.5,0) node[color=gray!50,above] {$\vect{x}_1$};
\draw[gray!50,->] (0,-2.5) -- (0,4.9) node[color=gray!50,above] {$\vect{x}_2$};

\node[coordinate] (g1) at (0,-2.5) {};
\node[coordinate] (g12) at (0,0) {};
\node[coordinate] (g2) at (2,2) {};
\node[coordinate] (g23) at (2,4.5) {};
\node[coordinate] (g34) at (-2,4.5) {};
\node[coordinate] (g4) at (-2,1) {};
\node[coordinate] (g41) at (0,1) {};
\node[coordinate] (g5) at (0,-3.5) {};

\filldraw[convex] (g34) -- (g4) -- (g41) -- (g1) -- (g12) -- (g2) -- (g23);
\draw[dashed,line width=1.8pt] (g5) -- (g1);
\draw[convexborder] (g34) -- (g4) -- (g41) -- (g1) -- (g12) -- (g2) -- (g23);
\filldraw[point] (g4) circle (3pt) node[left=3pt] {$\vect{p}^1$};
\filldraw[point] (g2) circle (3pt) node[right=3pt] {$\vect{p}^2$};
\filldraw[point] (g5) circle (3pt) node[left=3pt] {$\vect{p}^3$};
\draw[ray,line width=2pt,->,>=triangle 45] (0,2) -- node[right=3pt] {$\vect{r}^0$} (0,4) ;
\end{scope}
\begin{scope}[shift={(4.8,-2)},scale=0.8,>=triangle 45,convex/.style={draw=lightgray,fill=lightgray,fill opacity=0.7},point/.style={blue!50},line/.style={blue!50,ultra thick},convexborder/.style={ultra thick},pointlabel/.style={black},ray/.style={green!60!black}]
\equilateral{7}{90};

\barycenter{g1}{\expo{0}}{0}{\expo{0}};
\barycenter{g12}{\expo{0}}{\expo{0}}{\expo{0}};
\barycenter{g2}{\expo{2}}{\expo{2}}{\expo{0}};
\barycenter{g3}{0}{\expo{0}}{0};
\barycenter{g4}{\expo{-2}}{\expo{1}}{\expo{0}};
\barycenter{g41}{\expo{0}}{\expo{1}}{\expo{0}};

\filldraw[convex] (g1) -- (g12) -- (g2) -- (g3) -- (g4) -- (g41) -- cycle;
\draw[convexborder] (g1) -- (g12) -- (g2) -- (g3) -- (g4) -- (g41) -- cycle;
\filldraw[point] (g4) circle (3pt) node[above left=1pt] {$\vect{v}^1$};
\filldraw[point] (g2) circle (3pt) node[left=1pt] {$\vect{v}^2$};
\filldraw[point] (g1) circle (3pt) node[above right=1pt] {$\vect{v}^3$};
\filldraw[ray] (g3) circle (3pt) node[above right=1pt] {$\vect{v}^0$};
\end{scope}
\end{tikzpicture}
\end{center}
\caption{A tropical polyhedron in $\maxplus^2$ (left), and an equivalent representation by a cone in $\maxplus^3$ (right)}\label{fig:tropical_poly}
\end{figure}
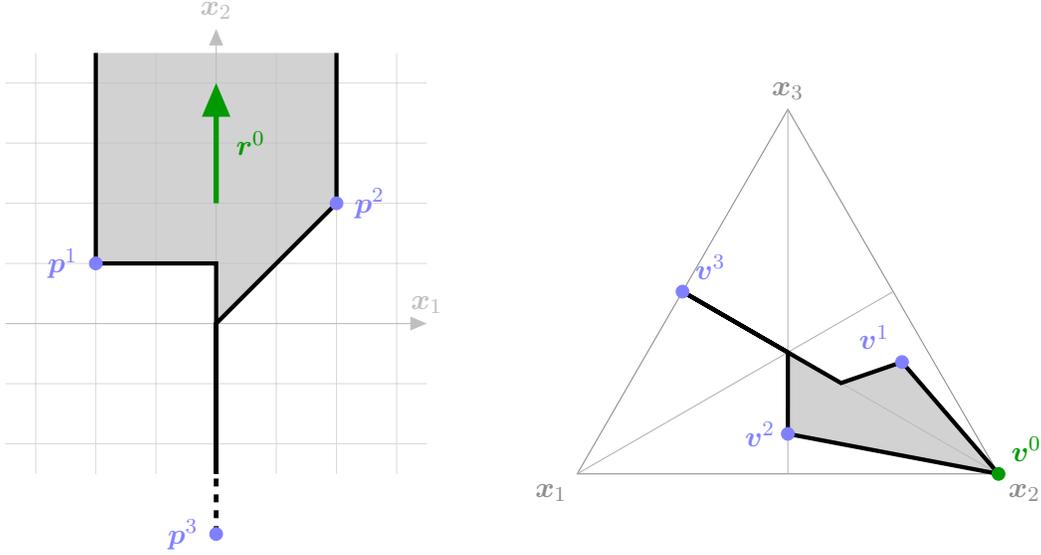

\begin{example}\label{ex:running}
In the following sections, we will illustrate our results on the tropical polyhedron $\PP$ depicted in solid gray (the black border is included) in the left hand side of Figure~\ref{fig:tropical_poly}. It is defined as the intersection of the half-spaces given by the inequalities:
\[
\begin{aligned}
0 & \leq \vect{x}_1+2 \\
\vect{x}_1 & \leq \max(\vect{x}_2,0) \\
\vect{x}_1 & \leq 2 \\
0 & \leq \max(\vect{x}_1,\vect{x}_2-1)
\end{aligned}
\]
This polyhedron is generated by the vertices $\vect{p}^1 = (-2,1)$, $\vect{p}^2 = (2,2)$, and $\vect{p}^3 = (0,-\infty)$, and by the extreme ray $\mpray{\vect{r}^0}$ where $\vect{r}^0 = (-\infty,0)$.
	
Its homogenized cone $\CC$ is depicted in the right hand side of Figure~\ref{fig:tropical_poly}. This cone is represented in barycentric coordinates: each element $(\vect{x}_1,\vect{x}_2,\vect{x}_3)$ is represented as a barycenter with weights $(e^{\vect{x}_1}, e^{\vect{x}_2}, e^{\vect{x}_3})$ of the three vertices of the outermost triangle. Two representatives of a same ray are thus represented by the same point. Besides, this barycentric representation is convenient to represent points with infinite coordinates, which are mapped to the boundary of the triangle. The cone $\CC$ is given by the linear inequalities:
\[
\begin{aligned}
\vect{x}_3 & \leq \vect{x}_1+2 \\
\vect{x}_1 & \leq \max(\vect{x}_2,\vect{x}_3) \\
\vect{x}_1 & \leq \vect{x}_3+2 \\
\vect{x}_3 & \leq \max(\vect{x}_1,\vect{x}_2-1)
\end{aligned}
\label{eq:running_cone}
\]
In accordance with Proposition~\ref{prop:homo_gen}, it is generated by the extreme elements $\vect{v}^0 = (-\infty,0,-\infty)$, $\vect{v}^1 = (-2,1,0)$, $\vect{v}^2 = (2,2,0)$, and $\vect{v}^3 = (0,-\infty,0)$.
\end{example}

\section{Combinatorial characterization of extremality using directed hypergraphs}\label{sec:combinatorial_characterization}

We first show that the extremality of a vector of a (non necessarily polyhedral) tropical cone is a local property.
\begin{propositiondefinition}\label{prop:local_extremality}
Let $\CC$ be a tropical cone, and $\vect{v}$ be a non-null vector of $\CC$. Then $\vect{v}$ is extreme in $\CC$ if, and only if, there exists a neighborhood $N$ of $\vect{v}$ such that:
\[
\forall \vect{x}, \vect{y} \in \CC \cap N, \ \vect{v} = \vect{x} \mpplus \vect{y} \implies \vect{v} = \vect{x} \text{ or } \vect{v} = \vect{y}.
\]
In the latter case, $\vect{v}$ is said to be \emph{locally extreme} in $\CC$.
\end{propositiondefinition}

\begin{proof}
The ``only if'' part of the result is straightforward. Suppose that $\vect{v}$ is locally extreme in $\CC$, and let $N$ be a neighborhood of $\vect{v}$ as in Proposition~\ref{prop:local_extremality}. Suppose that $\vect{x}, \vect{y} \in \CC$ are two vectors distinct from $\vect{v}$, and satisfying $\vect{v} = \vect{x} \mpplus \vect{y}$. Consider $\alpha > 0$ sufficiently small so that $\vect{x}' \defi \alpha^{-1} \vect{v} \mpplus \vect{x}$ and $\vect{y}' \defi \alpha^{-1} \vect{v} \mpplus \vect{y}$ both belong to $N$. Clearly, $\vect{v} = \vect{x}' \mpplus \vect{y}'$, and it follows that $\vect{v}$ is equal to one of the two vectors $\vect{x}'$ and $\vect{y}'$. As $\alpha > 0$, this yields a contradiction.
\end{proof}

The \emph{support} of a vector $\vect{x}=(\vect{x}_i) \in \maxplus^d$
is defined as the set of the indices of its non-null coordinates:
\[
\supp(\vect{x}) \defi \{ i \in \oneto{d} \mid \vect{x}_i \neq \mpzero \}.
\]
The following proposition states that the extremality of an element of a tropical cone can be established only by considering the vectors of the cone which have a smaller support:
\begin{proposition}\label{prop:support_extremality}
Let $\CC \subset \maxplus^d$ be a tropical cone, and $\vect{v}$ be a non-null vector of $\CC$. Then the following two statements are equivalent:
\begin{enumerate}[(i)]
\item \label{item:a} $\vect{v}$ is extreme in $\CC$,
\item \label{item:b} $\vect{v}$ is extreme in $\{ \vect{x} \in \CC \mid \supp(\vect{x}) \subset \supp(\vect{v})  \}$.
\end{enumerate}
\end{proposition}

\begin{proof}
Let $\DD \defi \{  \vect{x} \in \CC \mid \supp(\vect{x}) \subset \supp(\vect{v}) \}$. It is straightforward that $\DD$ is a tropical cone. Besides, for any vectors $\vect{x}, \vect{y} \in \CC$ such that $\vect{v} = \vect{x} \mpplus \vect{y}$, their supports are both included in $\supp(\vect{v})$, hence $\vect{x}, \vect{y} \in \DD$. This concludes the proof.
\end{proof}

For the rest of the section, the set $\CC$ is supposed to be a polyhedral cone defined by a system of linear inequalities $A \vect{x} \leq B \vect{x}$, with $A, B \in \maxplus^{p \times d}$. Using Proposition~\ref{prop:support_extremality}, it is assumed that $\vect{v}$ is an element of $\CC$ satisfying $\supp(\vect{v}) = \oneto{d}$, up to considering the extracted system $A' \vect{y} \leq B' \vect{y}$, where $A'$ and $B'$ are respectively the matrices formed by the columns of $A$ and $B$ of index $i$ in the support of $\vect{v}$. We denote by $A_k$ (resp.\ $B_k$) the $k$-th row of the matrix $A$ (resp.\ $B$), for $k \in [p]$. 

Following the line of Proposition~\ref{prop:local_extremality}, we introduce the notion of tangent cone at the point $\vect{v}$, which captures the constraints induced by the cone $\CC$ in the neighborhood of $\vect{v}$.
\begin{definition}\label{def:tangent_cone}
The \emph{tangent cone} to $\CC$ at the element $\vect{v}$ is the tropical polyhedral cone $\tangent(\vect{v},\CC)$ defined by the following intersection of half-spaces:
\[
\tangent(\vect{v},\CC) \defi \;\;\bigcap_{
\mathclap{\substack{k\in[p]\\
A_k \vect{v} = B_k \vect{v} > \mpzero}}} \;\;
\Bigl\{\, \vect{x} \in \maxplus^d \bigm| \; \; \; \smashoperator{\mpsum_{i \in \argmax(A_k \vect{v})}} \; \; \vect{x}_i \; \; \leq \; \; \; \; \smashoperator{\mpsum_{j \in \argmax(B_k \vect{v})}} \; \; \vect{x}_j \, \Bigr\},\label{eq:tangent_cone}
\]
where for any $\vect{c} = (\vect{c}_i) \in \maxplus^{1 \times d}$, the set $\argmax(\vect{c} \vect{v})$ is defined as the argument of the maximum $\vect{c} \vect{v} = \max_{i\in[d]} (\vect{c}_i + \vect{v}_i)$.
\end{definition}

\begin{proposition}\label{prop:tangent_cone}
There exists a neighborhood $N$ of $\vect{v}$ such that for all $\vect{x} \in N$, $\vect{x}$ belongs to $\CC$ if, and only if, it is an element of $\vect{v} + \tangent(\vect{v},\CC)$.
\end{proposition}

\begin{proof}
Consider a neighborhood $N$ in which all elements $\vect{x}$ satisfy the following conditions:
\begin{enumerate}[(i)]
\item $A_k \vect{x} < B_k \vect{x}$ for all $k \in [p]$ such that $A_k \vect{v} < B_k \vect{v}$, 
\item $\argmax(A_k \vect{x}) \subset \argmax(A_k \vect{v})$ and $\argmax(B_k \vect{x}) \subset \argmax(B_k \vect{v})$ for any other $k \in [p]$.
\end{enumerate}
Let $\vect{x} \in N$. Note that $\vect{x}$ belongs to $\CC$ if, and only if, for each $k \in [p]$ verifying $A_k \vect{v} = B_k \vect{v}$, 
\begin{equation}
\max_{i \in \argmax(A_k \vect{v})} (A_{ki} + \vect{x}_i) \leq \max_{j \in \argmax(B_k \vect{v})} (B_{kj} + \vect{x}_j), \label{eq:eq1}
\end{equation}
by definition of $N$.

Suppose that $\vect{x}$ belongs to $\CC$. Let $k$ such that $A_k \vect{v} = B_k \vect{v} > \mpzero$. Since for all $i \in \argmax(A_k \vect{v})$ and $j \in \argmax(B_k \vect{v})$, $A_{ki} + \vect{v}_i = B_{kj} + \vect{v}_j > \mpzero$, the term $A_{ki} + \vect{v}_i$ (resp.\ $B_{kj} + \vect{v}_j$) can be subtracted from $A_{ki} + \vect{x}_i$ (resp.\ $B_{kj} + \vect{x}_j$) in~\eqref{eq:eq1}, which shows:
\begin{equation}
\max_{i \in \argmax(A_k \vect{v})} (\vect{x}_i - \vect{v}_i) \leq \max_{j \in \argmax(B_k \vect{v})} (\vect{x}_j - \vect{v}_j). \label{eq:eq2}
\end{equation}

Conversely, suppose that $\vect{x} - \vect{v} = (\vect{x}_i - \vect{v}_i)_i$ is an element of $\tangent(\vect{v}, \CC)$. Consider $k \in [p]$ such that $A_k \vect{v} = B_k \vect{v} > \mpzero$. Adding the term $A_{ki} + \vect{v}_i$ (resp.\ $B_{kj} + \vect{v}_j$) to each $\vect{x}_i - \vect{v}_i$ (resp.\ $\vect{x}_j - \vect{v}_j$) in~\eqref{eq:eq2} shows that $\vect{x}$ satisfies~\eqref{eq:eq1}. Besides, if $A_k \vect{v} = B_k \vect{v} = \mpzero$, the row vectors $A_k$ and $B_k$ are identically null (since $\vect{v}$ has a full support), and~\eqref{eq:eq1} is trivially satisfied.
\end{proof}

\begin{remark}
The term \emph{tangent cone} is borrowed from convex analysis and optimization, where it is usually defined as the set of the directions which are asymptotically admissible from the vector $\vect{v}$ in the set $\CC$:
\[
\tangent^0(\vect{v},\CC) \defi \Bigl\{ \vect{x} \in \mathbb{R}^d \mid \exists \{t_k\}_{k \geq 0} \in (\mathbb{R}_+^*)^\mathbb{N}, \{\vect{v}^k\}_{k \geq 0} \in \CC^\mathbb{N}, \ t_k \mathop{\longrightarrow}\limits_{k\to +\infty} 0, \ (\vect{v}^k-\vect{v})/t_k \mathop{\longrightarrow}\limits_{k\to +\infty} \vect{x} \Bigr\}. \label{eq:usual_tangent_cone_def}
\]
We claim that $\tangent(\vect{v},\CC)$ defined here is the topological closure of $\tangent^0(\vect{v},\CC)$ (the latter is closed in $\mathbb{R}^d$, but not in $\maxplus^d$). Using essentially the same technique as in the proof of Proposition~\ref{prop:tangent_cone}, it can be shown that $\tangent^0(\vect{v},\CC)$ coincides with the set $\tangent(\vect{v},\CC) \cap \mathbb{R}^d$. Thus, the inclusion $\closure(\tangent^0(\vect{v},\CC)) \subset \tangent(\vect{v},\CC)$ is implied by the fact that $\tangent(\vect{v},\CC)$ is closed, as any tropical polyhedral cone. The opposite inclusion $\tangent(\vect{v},\CC) \subset \closure(\tangent^0(\vect{v},\CC))$ comes from  that the vector~$\mponevect$ necessarily belongs to $\tangent(\vect{v},\CC)$, hence any element $\vect{x} \in \tangent(\vect{v},\CC)$ can be expressed as the limit of the sequence of the elements $\vect{x}^k \defi \vect{x} \mpplus (k^{-1} \mponevect)$, for $k$ tending to $+\infty$. As for $k \geq 0$, every vector $\vect{x}^k$ belongs to $\tangent(\vect{v},\CC) \cap \mathbb{R}^d = \tangent^0(\vect{v},\CC)$, this completes the proof of the claim.
\end{remark}

Combining Proposition~\ref{prop:tangent_cone} with the local characterization of extremality given by Proposition~\ref{prop:local_extremality} yields the following reduction:
\begin{proposition}\label{prop:unit_extremality}
The vector $\vect{v}$ is extreme in $\CC$ if, and only if, the element $\mponevect$ is extreme in $\tangent(\vect{v}, \CC)$.
\end{proposition}

\begin{proof}
Let $N'$ be the set consisting of the elements $\vect{x} - \vect{v}$ for $\vect{x} \in N$, where $N$ is given by Proposition~\ref{prop:tangent_cone}. 
First remark that 
$\mponevect \in \tangent(\vect{v},\CC)$.

Observe that two vectors $\vect{x}, \vect{y} \in \CC \cap N$ satisfy $\vect{v} = \vect{x} \mpplus \vect{y}$ if, and only if, $\vect{x'} \defi \vect{x} - \vect{v}$ and $\vect{y'} \defi \vect{y} - \vect{v}$ belongs to $\tangent(\vect{v},\CC) \cap N'$ (by Proposition~\ref{prop:tangent_cone}), and $\mponevect = \vect{x}' \mpplus \vect{y}'$. We deduce that $\vect{v}$ is locally extreme in $\CC$ if, and only if, $\mponevect$ is locally extreme in $\tangent(\vect{v},\CC)$. We conclude the proof by Proposition~\ref{prop:local_extremality}.
\end{proof}
 
The interest of Proposition~\ref{prop:unit_extremality} is that we are now reduced to characterizing extremality of a vector of $\BB^d$ in a closed tropical cone $\DD$ (here $\tangent(\vect{v},\CC)$), which is stable under the usual multiplication by positive scalars, \ie{} $\alpha \times \vect{x} \in \DD$ for all $\vect{x} \in \DD$ and $\alpha > 0$ (the multiplication being understood entrywise). 

\begin{proposition}\label{prop:usualcone_extremality}
Let $\DD \subset \maxplus^d$ be a closed tropical cone, which is stable under the multiplication by positive scalars in the usual sense. Then a vector $\vect{w} \in \DD \cap \BB^d$ is extreme in $\DD$ if, and only if, it is extreme in set $\DD \cap \BB^d$ seen as a boolean cone over $\BB^d$.
\end{proposition}

\begin{proof}
The ``only if''	part is obvious. 

For the ``if'' part, suppose that $\vect{w} = \vect{x} \mpplus \vect{y}$, where $\vect{x}, \vect{y} \in \DD$. Observe that every entry of $\vect{x}$ and $\vect{y}$ is non-positive, and consider the sequences of the vectors $\vect{x}^k \defi k \times \vect{x}$ and $\vect{y}^k \defi k \times \vect{y}$, for $k \geq 1$. By assumption, all elements $\vect{x}^k$ and $\vect{y}^k$ belong to $\DD$. Besides, $\vect{x}^k \mpplus \vect{y}^k = k \times \vect{w} = \vect{w}$, since $\vect{w} \in \BB^d$. The two sequences $(\vect{x}^k)_k$ and $(\vect{y}^k)_k$ admit a limit, respectively denoted by $\vect{x}'$ and $\vect{y}'$, which both belong to $\DD \cap \BB^d$ (since $\DD$ is topologically closed). The element $\vect{w}$ being extreme in the boolean cone $\DD \cap \BB^d$, one of the two vectors $\vect{x}'$ or $\vect{y}'$ is equal to $\vect{w}$. Supposing for instance that $\vect{w} = \vect{x}'$, we know that $\vect{x}' \leq \vect{x} \leq \vect{w}$, and we conclude that $\vect{w} = \vect{x}$.
\end{proof}

Instantiating Proposition~\ref{prop:usualcone_extremality} with the tangent cone $\tangent(\vect{v},\CC)$ provides the following combinatorial characterization of the extremality of $\vect{v}$:
\begin{theorem}\label{th:tangent_cone}
Let $\CC \subset \maxplus^d$ be a polyhedral cone, and $\vect{v}$ a vector of $\CC$ with full support. The following three propositions are equivalent:
\begin{enumerate}[(i)]
\item\label{item:i} the vector $\vect{v}$ is extreme in $\CC$,
\item\label{item:ii} the vector $\mponevect$ is extreme in the boolean cone $\tangent(\vect{v},\CC) \cap \BB^d$,
\item\label{item:iii} there exists $i \in \oneto{d}$ such that the following inequalities hold:
\begin{equation}
\forall j \in [d], \ \forall \vect{x} \in \tangent(\vect{v},\CC) \cap \BB^d, \ \vect{x}_i \leq \vect{x}_j.  \label{eq:zero_propagation}
\end{equation}
\end{enumerate}
\end{theorem}

\begin{proof}
The equivalence $\eqref{item:i} \Leftrightarrow \eqref{item:ii}$ follows from Propositions~\ref{prop:unit_extremality} and~\ref{prop:usualcone_extremality}. It remains to show $\eqref{item:ii} \Leftrightarrow \eqref{item:iii}$.

Suppose that $\mponevect$ is extreme in the boolean cone $\DD \defi \tangent(\vect{v}, \CC) \cap \BB^d$. We claim that there exists $i \in [d]$ such that for all $\vect{x} \in \DD \setminus \{\mponevect\}$, $\vect{x}_i = \mpzero$. If not, $\mponevect$ could be written as the sum of some elements $\vect{x}^j \in \DD$ such that $\vect{x}^j_j = \mpone$, for $j = 1, \dots, d$. It follows that for all elements $\vect{x}$ of $\DD$ and $j \in [d]$, $\vect{x}_i \leq \vect{x}_j$.

Reciprocally, if~\eqref{eq:zero_propagation} is satisfied, then every vector $\vect{x} \in \DD \setminus \{\mponevect\}$ verifies $\vect{x}_i = \mpzero$. Then for all pairs $(\vect{x}, \vect{y})$ of such elements, we have $(\vect{x} \mpplus \vect{y})_i = \mpzero$. We conclude that $\mponevect$ is extreme.
\end{proof}

Since the tangent cone $\tangent(\vect{v},\CC)$ is defined by inequalities with coefficients in $\BB$, the nature of the characterization provided by Theorem~\ref{th:tangent_cone} is purely boolean. However, testing, by exploration, whether there exists $i \in [d]$ such that every vector of $\tangent(\vect{v},\CC) \cap \BB^d$ satisfies $\vect{x}_i \leq \vect{x}_j$ for all $j \in [d]$, does not have acceptable complexity. Instead, we propose to express the satisfiability of the inequalities~\eqref{eq:zero_propagation} as a reachability problem on directed hypergraphs.

Directed hypergraphs are a generalization of directed graphs, in which arcs leave and enter subsets of vertices. A \emph{directed hypergraph} over the nodes $1, \dots, d$ is a set of \emph{hyperarcs} of the form $(T, H)$, where $T, H \in 2^d$. The notion of reachability is extended from directed graphs to directed hypergraphs, and defined inductively as follows: given a hypergraph $\GG$, and $i, j \in [d]$, $i$ is reachable from $j$ in $\GG$, which is denoted by $j \reach_\GG i$, if $i = j$, or there exists a hyperarc $(T, H)$ in $\GG$ such that $i \in H$, and all elements of $T$ are reachable from $j$. 

\begin{figure}[t]
\begin{minipage}[t]{0.49\textwidth}
\begin{center}
\begin{tikzpicture}[>=triangle 45,scale=0.8, vertex/.style={circle,draw=black,very thick,minimum size=2ex}, hyperedge/.style={draw=black,thick}, simpleedge/.style={draw=black,thick}]
\node [vertex] (u) at (-2,-1) {$1$};
\node [vertex] (v) at (0,0) {$2$};
\node [vertex] (w) at (0,-2)  {$3$};
\node [vertex] (x) at (3.5,0) {$4$};
\node [vertex] (y) at (3.5,-2) {$5$};
\node [vertex] (t) at (2,-4.5) {$6$};

\path[->] (u) edge[simpleedge,out=90,in=-180] (v);
\node at (-1,-0.5) {$a_1$};
\path[->] (v) edge[simpleedge,out=-90,in=90] (w);
\node at (-0.5,-1.3) {$a_2$};
\path[->] (w) edge[simpleedge,out=-120,in=-60] (u);
\node at (-1.5,-2.5) {$a_3$};
\node at (1.75,-0.5) {$a_4$};
\node at (2.5,-3.5) {$a_5$};
\oldhyperedge{v}{v,w}{w}{0.5}{0.6}{y,x};
\oldhyperedge{y}{y,w}{w}{0.6}{0.6}{t};
\end{tikzpicture}
\end{center}
\caption{A directed hypergraph}\label{fig:hypergraph}
\end{minipage}
\begin{minipage}[t]{0.49\textwidth}
\begin{center}
\begin{tikzpicture}[scale=0.8,>=triangle 45,convex/.style={draw=lightgray,fill=lightgray},point/.style={blue!50},line/.style={blue!50,ultra thick},convexborder/.style={ultra thick},pointlabel/.style={black},ray/.style={ultra thick,green!60!black}]
\equilateral{7}{90};

\barycenter{g1}{\expo{0}}{0}{\expo{0}};
\barycenter{g12}{\expo{0}}{\expo{0}}{\expo{0}};
\barycenter{g2}{\expo{2}}{\expo{2}}{\expo{0}};
\barycenter{g3}{0}{\expo{0}}{0};
\barycenter{g4}{\expo{-2}}{\expo{1}}{\expo{0}};
\barycenter{g41}{\expo{0}}{\expo{1}}{\expo{0}};

\barycenter{g2y}{exp(1)}{0}{exp(0)};
\barycenter{g2z}{exp(1)}{exp(1)}{0};

\filldraw[convex] (g1) -- (g12) -- (g2) -- (g3) -- (g4) -- (g41) -- cycle;
\draw[convexborder,miter limit=1] (g1) -- (g12) -- (g2) -- (g3) -- (g4) -- (g41) -- cycle;

\filldraw[blue!30,draw=none,fill opacity=0.3] (z) -- (g2) -- (y);
\draw[blue!60, dotted, ultra thick] (z) -- (g2) -- (y) -- (z);

\node[circle,fill=black,inner sep=0pt,minimum size=1.25ex] at (g2) [label=left:$\vect{v}^2$] {};
\end{tikzpicture}
\end{center}
\caption{The set $\vect{v}^2+\tangent(\vect{v}^2,\CC)$ (in light blue)}\label{fig:tangent_cone}
\end{minipage}
\end{figure}

\begin{example}
Figure~\ref{fig:hypergraph} depicts an example of directed hypergraph consisting of the hyperarcs $a_1 =(\{ 1\}, \{2\})$, $a_2 = (\{ 2\}, \{3\})$, $a_3 = (\{ 3\}, \{1\})$, $a_4 = (\{ 2,3 \}, \{4,5\})$, and $a_5 = (\{ 3,5 \}, \{6\})$. 
We visualize a hyperarc as a bundle of arrows: a solid disk sector
indicates that the different arrows going through it belong to the same hyperarc; the head (resp.\ tail) of the hyperarc is the union of the 
heads (resp.\ tails) of these arrows. Applying the recursive definition of reachability from the node $1$ discovers the nodes $2$, then $3$, which leads to the two nodes $4$ and $5$ through the hyperarc $a_4$, and finally the node $6$ through $a_5$.
\end{example}

We now introduce the notion of \emph{tangent directed hypergraph}, which is an equivalent encoding of the tangent cone as a directed hypergraph. It derives from the system of inequalities defining the tangent cone:
\begin{definition}
The \emph{tangent directed hypergraph} at $\vect{v}$, denoted by $\GG(\vect{v},\CC)$, is the directed hypergraph consisting of the hyperarcs $(\argmax(B_k \vect{v}), \argmax(A_k \vect{v}))$ for every $k \in [p]$ such that $A_k \vect{v} = B_k \vect{v} > \mpzero$.
\end{definition}

The reachability relation of the tangent directed hypergraph precisely captures the constraints of the form $\vect{x}_i \leq \vect{x}_j$ satisfied by the boolean elements $\vect{x}$ of the tangent cone, as shown below:
\begin{proposition}\label{prop:hypergraph_reachability}
The following two conditions are equivalent:
\begin{enumerate}[(a)]
\item\label{item:alpha} $i$ is reachable from $j$ in the hypergraph $\GG(\vect{v},\CC)$,
\item\label{item:beta} for all $\vect{x} \in \tangent(\vect{v},\CC) \cap \BB^d$, $\vect{x}_i \leq \vect{x}_j$,
\end{enumerate}
\end{proposition}

\begin{proof}
We first prove by induction that for every node $i$ reachable from $j$, the inequality $\vect{x}_i \leq \vect{x}_j$ holds for all $\vect{x} \in \tangent(\vect{v},\CC) \cap \BB^d$. The case $i = j$ is trivial. Suppose that there exists a hyperarc $(T,H)$ in $\GG(\vect{v},\CC)$ such that $i \in H$, and each $k \in T$ is reachable from $j$. Then for all $\vect{x} \in \tangent(\vect{v},\CC) \cap \BB^d$, $\mpplus_{k \in T} \vect{x}_k \leq \vect{x}_j$. Besides, $\vect{x}_i \leq \mpplus_{k \in T} \vect{x}_k$ by definition of $\tangent(\vect{v},\CC)$ and $\GG(\vect{v},\CC)$. Hence $\vect{x}_i \leq \vect{x}_j$. This shows that $\eqref{item:alpha} \Rightarrow \eqref{item:beta}$.

Then, let us show by contraposition that $\eqref{item:beta} \Rightarrow \eqref{item:alpha}$. Let $R$ be the set of nodes $k \in [d]$ reachable from $j$, and assume that $i \not \in R$. Consider the element $\vect{x} \in \BB^d$ defined by $\vect{x}_k = \mpzero$ if $k \in R$, and $\mpone$ otherwise. In particular, $\vect{x}_i > \vect{x}_j$. Let $(T,H)$ be a hyperarc of $\GG(\vect{v},\CC)$. If $T \subset R$, then $H$ is included into $R$, so that the inequality $\mpplus_{k \in H} \vect{x}_k \leq \mpplus_{l \in T} \vect{x}_l$ is satisfied. If $T \not \subset R$, the latter inequality is still valid because $\mpplus_{l \in T} \vect{x}_l = \mpone$. We deduce that $\vect{x}$ belongs to $\tangent(\vect{v},\CC) \cap \BB^d$.
\end{proof}

\begin{remark}\label{remark:hypergraph_reachability}
The two statements of Proposition~\ref{prop:hypergraph_reachability} can also be shown to be equivalent to the property:
\begin{inparaenum}[(a)]
\setcounter{enumi}{2}
\item\label{eq:gamma} for all $\vect{x} \in \tangent(\vect{v},\CC)$, $\vect{x}_i \leq \vect{x}_j$.
\end{inparaenum}
Replacing $\tangent(\vect{v},\CC) \cap \BB^d$ by $\tangent(\vect{v},\CC)$, the first part of the proof of Proposition~\ref{prop:hypergraph_reachability} indeed shows $\eqref{item:alpha} \Rightarrow \eqref{eq:gamma}$. The implication $\eqref{eq:gamma} \Rightarrow \eqref{item:beta}$ is trivial. 
\end{remark}

Given a directed hypergraph $\GG$, the \emph{strongly connected components} (\scc{}s for short) are defined as the equivalence classes of the relation $\reacheq_\GG$, given by $i \reacheq_\GG j$ if $i \reach_\GG j$ and $j \reach_\GG i$. Strongly connected components are partially ordered by the relation $\reachleq_\GG$ induced by reachability, \ie{} $C \reachleq_\GG C'$ if $C$ and $C'$ admit a representative $i$ and $j$ respectively such that $i \reach_\GG j$. 

The theorem stated in the introduction now follows as a consequence of the previous results.
\begin{proof}[\proofname{} of Theorem~\ref{th-main}]
From Theorem~\ref{th:tangent_cone} and Proposition~\ref{prop:hypergraph_reachability}, the vector $\vect{v}$ is extreme if, and only if, there is a node $i$ reachable from every node $j$ in the tangent directed hypergraph $\GG(\vect{v},\CC)$. This holds if, and only if, $i$ belongs to a strongly connected component $C$ such that $C' \reachleq_{\GG(\vect{v},\CC)} C$ for any \scc{} $C'$.
\end{proof}

\begin{example}\label{ex:running_tangent}
Let us illustrate Theorem~\ref{th:main} by establishing the extremality of the vector $\vect{v}^2 = (2,2,0)$ of the cone $\CC$ defined in Example~\ref{ex:running}. In~\eqref{eq:running_tangent1}, the inequalities which are active at $\vect{v}^2$ are colored in red, and the terms which belong to the arguments of the left-/right-hand sides are underlined. 
A system of inequalities defining the cone $\tangent(\vect{v}^2,\CC)$, in~\eqref{eq:running_tangent2}, is obtained by keeping only the underlined terms.
\begin{center}
\vskip-1.5ex
\begin{minipage}[c]{0.3\textwidth}
\vspace{0ex}
\begin{equation}
\begin{aligned}
\vect{x}_3 & \leq \vect{x}_1+2 \\
\textcolor{red!80!black}{\underline{\vect{x}_1}} & \textcolor{red!80!black}{{}\leq \max(\underline{\vect{x}_2},\vect{x}_3)} \\
\textcolor{red!80!black}{\underline{\vect{x}_1}} & \textcolor{red!80!black}{{}\leq \underline{\vect{x}_3+2}} \\
\vect{x}_3 & \leq \max(\vect{x}_1,\vect{x}_2-1)
\end{aligned}
\label{eq:running_tangent1}
\end{equation}
\end{minipage}
\qquad
\begin{minipage}[c]{0.2\textwidth}
\vspace{0ex}
\begin{equation}
\begin{aligned}
\vect{x}_1 & \leq \vect{x}_2 \\
\vect{x}_1 & \leq \vect{x}_3 
\end{aligned}
\label{eq:running_tangent2}
\end{equation}
\end{minipage}
\end{center}
\vskip1ex
Figure~\ref{fig:tangent_cone} illustrates that the cones $\CC$ and $\vect{v}^2 + \tangent(\vect{v}^2,\CC)$ locally coincide in a neighborhood of $\vect{v}^2$.	
The tangent directed hypergraph $\GG(\vect{v}^2,\CC)$ associated with the vector $\vect{v}^2$ is formed by the two hyperarcs $(\{2\},\{1\})$ and $(\{3\},\{1\})$. The node $1$ consequently forms the greatest strongly connected component of the hypergraph (for the partial order $\reachleq_{\GG(\vect{v}^2,\CC)}$).
\end{example}

\begin{remark}\label{remark:minimality}
It was shown in the proof of Theorem~3.1 of~\cite{GK} (see also~\cite{GK06a}), and independently in~\cite[Theorem~14]{BSS}, that the vector $\vect{v} \in \CC$ is extreme if, and only, if there exists $i \in \oneto{d}$ such that $\vect{v}$ is \emph{minimal of type $i$}, \ie{} minimal in the set $\{ \vect{x} \in \CC \mid \vect{x}_i = \vect{v}_i\}$. This result can be recovered as a corollary of Theorem~\ref{th:main}.

Observe that $i$ is reachable from any node $j \in [d]$ in $\GG(\vect{v},\CC)$ if, and only if, the hypergraph $\GG' = \GG(\vect{v},\CC) \cup \{\, (\{i\}, \{j\}) \mid j \in [d] \, \}$ is strongly connected. Let $\CC' = \{\vect{x} \in \maxplus^d \mid A \vect{x} \leq B \vect{x}, \ \vect{x}_j - \vect{v}_j \leq \vect{x}_i - \vect{v}_i \text{ for all } j \in [d]\}$. Then $\GG'$ is precisely the tangent directed hypergraph $\GG(\vect{v}, \CC')$. By Proposition~\ref{prop:hypergraph_reachability} and Remark~\ref{remark:hypergraph_reachability}, $\GG'$ is strongly connected if, and only if, the tangent cone $\tangent(\vect{v},\CC')$ is reduced to the ray $\maxplus \mponevect$. By Proposition~\ref{prop:tangent_cone}, this amounts to the equality $\CC' \cap N = \maxplus \vect{v}$ for a certain neighborhood $N$ of $\vect{v}$, or equivalently, $\CC' = \maxplus \vect{v}$ (as $\CC'$ is connected). The latter holds if, and only if, $\vect{v}$ is minimal of type $i$.

It follows that the nodes contained in the greatest strongly connected component of $\GG(\vect{v},\CC)$ (if it exists) are precisely the integers $i \in \oneto{d}$ such that $\vect{v}$ is minimal of type $i$. This suggests an alternative proof of Theorem~\ref{th:main}.
\end{remark}

An algorithm due to Gallo \etal{}~\cite{GalloDAM93} shows that the set of nodes that are reachable from a given node in a directed hypergraph $\GG$ can be computed in linear time in the size of the hypergraph, $\size(\GG) \defi \sum_{(T,H) \in \GG} (\card{T} + \card{H})$. The following result shows that one can in fact compute the maximal \scc{}s  with almost the same complexity. The algorithm, which is too technical to be included here, is detailed in~\cite{AllamigeonAlgorithmica2013}. 
\begin{theorem}[\cite{AllamigeonAlgorithmica2013}]\label{th:complexity}
The set of maximal \scc{}s of a hypergraph $\GG$ over the nodes $1, \dots, d$ can be computed in time $O(\size(\GG) \times \alpha(d))$, where $\alpha$ denotes the inverse of the Ackermann function.
\end{theorem}
The function $\alpha$ is a very slowly growing map. In particular, $\alpha(x) \leq 5$ for any practical values of $x$. Hence, the complexity is said to be \emph{almost linear}. The term $\alpha(d)$ originates from the use of Tarjan's union-find structure~\cite{Tarjan84} to efficiently manipulate partitions of the set of nodes $1, \dots, d$. In the sequel, \Call{MaxScc}{} denotes an algorithm returning the set of the maximal \scc{}s of the hypergraph given in input, with the time complexity given in Theorem~\ref{th:complexity}.

As a consequence, the criterion of Theorem~\ref{th:main} can be very efficiently evaluated, in almost linear time in the size of the system of inequalities defining the tropical polyhedral cone. It can also benefit from the sparsity of the system, since the size of the tangent directed hypergraph is bounded by the number of non-null coefficients in the inequalities.

\section{The tropical double description method}\label{sec:external_to_internal}

We next present a tropical analogue of the double description method of Motzkin \etal{}~\cite{MRTT53}. The tropical method
computes a minimal generating set of a polyhedral cone, starting from a system of tropically linear inequalities defining it. We first deal with the inductive scheme of the method (Section~\ref{subsec:ddm}), then present the main algorithm (Section~\ref{subsec:resulting_algo}), and its extension to intersections with tropical hyperplanes (Section~\ref{subsec:hyperplanes}).

\subsection{Inductive scheme}\label{subsec:ddm}

The tropical double description method relies on an incremental technique based on a successive elimination of inequalities. Given a polyhedral cone defined by a system of $p$ constraints, it computes by induction on $k = 1, \dots, p$ a generating set $V_k$ of the intermediate cone defined by the first $k$ constraints. 

Passing from the set $V_k$ to the set $V_{k+1}$ relies on a result which, given a polyhedral cone $\CC$ and a tropical half-space $\HH$, allows to build a generating set $V'$ of $\CC \cap \HH$ from a generating set $V$ of $\CC$. This is referred to as the \emph{elementary step} of the scheme. 
(Note that the next result applies more generally to {\em non polyhedral}
tropical cones provided they are closed.)
\begin{theorem}[Elementary step]\label{th:ddm_basis}
Let $\CC \subset \maxplus^d$ be a closed tropical cone generated by a set $V$ of elements of $\maxplus^d$, and let $\HH$ be a half-space. Then the cone $\CC \cap \HH$ is generated by the following set:
\begin{equation}
(V \cap \HH) \cup \bigl\{  \vect{v} \mpplus \rho \vect{w} \mid (\vect{v},\vect{w})  \in (V \cap \HH) \times (V \setminus \HH), \ \rho = \max \{ \mu \lambda^{-1} \mid \lambda \vect{v} \mpplus \mu \vect{w} \in \HH \} \bigr\} . \label{eq:ddm_basis}
\end{equation}
\end{theorem}
The cone $\CC \cap \HH$ is thus generated by the elements $\vect{v}$ of $V$ satisfying the constraint associated with the half-space $\HH$, and by their pairwise combinations with the vectors $\vect{w}$ which are not located in $\HH$. Each combination $\vect{v} \mpplus \rho \vect{w}$ corresponds to the last element belonging to the half-space $\HH$ encountered when following the path from $\vect{v}$ to $\vect{w}$ along the tropical (projective) segment $\{ \lambda \vect{v} \mpplus \mu \vect{w} \mid \lambda, \mu \in \maxplus \}$.

\begin{remark}
Observe that the scalar $\rho$ is correctly defined in~\eqref{eq:ddm_basis}, meaning that the set $M = \{ \mu \lambda^{-1} \mid \lambda \vect{v} \mpplus \mu \vect{w} \in \HH\}$ admits a maximal element. First, note that $M$ is not empty (consider $\lambda = \mpone$ and $\mu = \mpzero$). We claim that it is bounded. If not, $\vect{v} \mpplus \rho \vect{w} \in \HH$ for arbitrarily large scalar $\rho \in \maxplus$, so that $\rho^{-1} \vect{v} \mpplus \vect{w} \in \HH$ as soon as $\rho > \mpzero$. Since $\HH$ is closed, this would imply that $\vect{w} \in \HH$ (taking $\rho \rightarrow +\infty$), which is nonsense. Finally, the set $M$ is closed, as the inverse image of $\HH$ by the continuous map $\rho \mapsto \vect{v} \mpplus \rho \vect{w}$. 
The supremum of the set $M$ therefore belongs to $M$.
\end{remark}

\begin{proof}[\proofname{} (Theorem~\ref{th:ddm_basis})]
Any element of the set given in~\eqref{eq:ddm_basis} obviously belongs to $\CC \cap \HH$.

Now consider $\vect{x} \in \CC \cap \HH$. Using the tropical analogue of the Minkowski-Carath{\'e}odory theorem established in~\cite[Theorem 3.1]{GK}, $\vect{x}$ can be written as a combination of at most $d$ elements of $V$, \ie{} there exist $V' \subset V \cap \HH$ and $W' \subset V \setminus \HH$ with $\card{V'} + \card{W'} \leq d$, and:
\[
\vect{x} = \mpsum_{\vect{v} \in V'} \lambda_{\vect{v}} \vect{v} \mpplus \mpsum_{\vect{w} \in W'} \lambda_{\vect{w}} \vect{w},
\]
where the $\lambda_{\vect{v}}$ and $\lambda_{\vect{w}}$ are non-zero scalars (in the tropical sense). Let $\rho_{(\vect{v},\vect{w})} = \max \{ \mu \lambda^{-1} \mid \lambda \vect{v} \mpplus \mu \vect{w} \in \HH \}$ for any pair $(\vect{v}, \vect{w}) \in V' \times W'$. First, let us show that for all $\vect{w} \in W'$, there exists $\vect{v} \in V'$ such that $\lambda_{\vect{v}} \rho_{(\vect{v},\vect{w})} \geq \lambda_{\vect{w}}$. If not, there is a given $\vect{w} \in W'$ satisfying $\lambda_{\vect{v}} \rho_{(\vect{v},\vect{w})} < \lambda_{\vect{w}}$ for all $\vect{v} \in V'$, hence $\lambda_{\vect{v}} \vect{v} \mpplus \lambda_{\vect{w}} \vect{w} \not \in \HH$ (since $\lambda_{\vect{v}} > \mpzero$) . Observe that the complementary of $\HH$ is stable by addition and by multiplication by a non-zero scalar, so that 
\[
\vect{x} = \Bigl(\mpsum_{\vect{v} \in V'} \lambda_{\vect{v}} \vect{v} \mpplus \lambda_{\vect{w}} \vect{w} \Bigr) \mpplus \Bigl(\mpsum_{\vect{w}' \in W' \setminus \{\vect{w}\}} \lambda_{\vect{w}'} \vect{w}'\Bigr)
\]
should not belong to $\HH$, which is a contradiction.

For all $\vect{w} \in W'$, let $\vect{v}_\vect{w}$ be an element of $V'$ such that $\lambda_{\vect{v}_\vect{w}} \rho_{(\vect{v}_\vect{w},\vect{w})} \geq \lambda_{\vect{w}}$. Since $\lambda_{\vect{w}}$ is not null, we have $\rho_{(\vect{v}_\vect{w},\vect{w})} > \mpzero$, hence we can write:
\[
\vect{x} = \Bigl(\mpsum_{\vect{v} \in V'} \lambda_{\vect{v}} \vect{v}\Bigr) \mpplus \Bigl(\mpsum_{\vect{w} \in W'} \bigl(\lambda_{\vect{w}} (\rho_{(\vect{v}_\vect{w},\vect{w})})^{-1}\bigr) (\vect{v}_\vect{w} \mpplus \rho_{(\vect{v}_\vect{w},\vect{w})} \vect{w})\Bigr).
\]
This shows that any element $\vect{x}$ of $\CC \cap \HH$ can be expressed as a combination of vectors of the set given in~\eqref{eq:ddm_basis}.
\end{proof}

Let us denote by $\vect{\epsilon}^i$ the element of $\maxplus^d$ whose $i$-th coordinate is equal to $\mpone$, and the other coordinates to $\mpzero$. The following theorem describes the whole inductive approach:
\begin{theorem}[Inductive scheme of the tropical double description method]\label{th:ddm}
Let $\CC \subset \maxplus^d$ be a polyhedral cone defined as the set $\{ \vect{x} \in \maxplus^d \mid A \vect{x} \leq B \vect{x}  \}$, where $A, B \in \maxplus^{p \times d}$ (with $p \geq 0$). Let $V_0,\dots,V_p$ be the sequence of finite subsets of $\maxplus^d$ defined as follows:
\[
\left\{
\begin{aligned}
V_0 & = \{\vect{\epsilon}^i\}_{i\in[d]},\\
V_k & = \bigl\{\vect{v} \in V_{k-1} \mid A_k \vect{v} \leq B_k \vect{v}\bigr\} \\
& \phantom{{}={}} {}\cup \bigl\{(A_k \vect{w}) \vect{v} \mpplus (B_k \vect{v}) \vect{w} \mid \vect{v},\vect{w} \in V_{k-1} \text{, } A_k \vect{v} \leq B_k \vect{v} \text{, and } A_k \vect{w} > B_k \vect{w} \bigr\} ,
\end{aligned}
\right.
\]
for all $k = 1, \dots, p$. 
Then $\CC$ is generated by the finite set $V_p$.
\end{theorem}

\begin{proof}
We show by using Theorem~\ref{th:ddm_basis} that each $V_k$ forms a generating set of the polyhedral cone $\{\vect{x} \in \maxplus^d \mid A_l \vect{x} \leq B_l \vect{x} \text{ for all } l = 1,\dots,k\}$. For $k = 0$, this is obvious.

Now suppose $k \geq 1$. Let $\HH_k$ be half-space defined by the inequality $A_k \vect{x} \leq B_k \vect{x}$, and $\vect{v}, \vect{w} \in V_{k-1}$ such that $\vect{v} \in \HH_k$ and $\vect{w} \not \in \HH_k$. We are going to show that $\max \{ \mu \lambda^{-1} \mid \lambda \vect{v} \mpplus \mu \vect{w} \in \HH_k \} = (A_k \vect{w})^{-1} (B_k \vect{v})$ (note that $A_k \vect{w}$ is not null as $A_k \vect{w} > B_k \vect{w}$). Indeed, if we set $\vect{x}$ to $\vect{v} \mpplus ((A_k \vect{w})^{-1} (B_k \vect{v})) \vect{w}$, we have $A_k \vect{x} = A_k \vect{v} \mpplus B_k \vect{v} \leq B_k \vect{v} \leq B_k \vect{x}$. Besides, if $\lambda > \mpzero$ and $\lambda \vect{v} \mpplus \mu \vect{w} \in \HH_k$, then 
$\lambda (A_k \vect{v}) \mpplus \mu (A_k \vect{w}) \leq \lambda (B_k \vect{v}) \mpplus \mu (B_k \vect{w})$. If $\mu > \mpzero$, then $\mu (A_k \vect{w}) > \mu (B_k \vect{w})$ so that $\mu (A_k \vect{w}) \leq \lambda (B_k \vect{v})$. Thus $\mu \lambda^{-1} \leq (A_k \vect{w})^{-1} (B_k \vect{w})$. The inequality also trivially holds as soon as $\mu = \mpzero$.

As a consequence, up to multiplicative factors, $V_k$ coincides with the set provided by Theorem~\ref{th:ddm_basis} for $V = V_{k-1}$ and $\HH = \HH_k$. This completes the proof.
\end{proof}

Observe that Theorem~\ref{th:ddm} provides a constructive proof of the ``Minkowski part'' of the Minkowski-Weyl theorem (Theorem~\ref{th:minkowski_weyl}), since it shows that all tropical polyhedral cones are generated by finite sets of elements of $\maxplus^d$. 

\begin{figure}
\begin{minipage}[t]{0.495\textwidth}
\vspace{0ex}
\begin{center}
\begin{tikzpicture}
[scale=0.8,>=triangle 45,
green/.style={circle,text=green!60!black,fill=green!60!black,inner sep=0pt,minimum size=1.25ex},
vtxlbl/.style={label distance=0pt,green!60!black},
red/.style={circle,text=red!80!black,fill=orange!30!red!90!black,inner sep=0pt,minimum size=1.25ex},
raylbl/.style={label distance=0pt,orange!30!red!90!black},
new/.style={circle,text=black,fill=blue!70,inner sep=0pt,minimum size=1.25ex},
newlbl/.style={label distance=0pt,blue!70}]
\equilateral{7}{90};

\barycenter{g1}{\expo{0}}{0}{\expo{0}};
\barycenter{g12}{\expo{0}}{\expo{0}}{\expo{0}};
\barycenter{g2}{\expo{2}}{\expo{2}}{\expo{0}};
\barycenter{g3}{0}{\expo{0}}{0};
\barycenter{g4}{\expo{-2}}{\expo{1}}{\expo{0}};
\barycenter{g41}{\expo{0}}{\expo{1}}{\expo{0}};

\barycenter{c1}{\expo{0}}{0}{0};
\barycenter{c2}{0}{\expo{2.5}}{\expo{0}};
\barycenter{c3}{0}{0}{\expo{0}};

\filldraw[lightgray] (g1) -- (g12) -- (g2) -- (g3) -- (g4) -- (g41) -- cycle;
\draw[draw=black,ultra thick] (g1) -- (g12) -- (g2) -- (g3) -- (g4) -- (g41) -- cycle;

\filldraw[blue!20,draw=none,fill opacity=0.3] (c1) -- (c2) -- (c3) -- cycle;
\draw[blue!40, ultra thick] (c1) -- (c2);

\node[red] at (g3) [label={[raylbl]above:$\vect{v}^0$}] {};
\node[green] at (g4) [label={[vtxlbl,label distance=-0.5ex]above left:$\vect{v}^1$}] {};
\node[green] at (g2) [label={[vtxlbl]left:$\vect{v}^2$}] {};
\node[green] at (g1) [label={[vtxlbl]above right:$\vect{v}^3$}] {};

\node[new] at (intersection cs:
  first line={(c1) -- (c2)},
  second line={(g3) -- (g2)}) [label={[newlbl,label distance=-0.5ex]below:$\vect{w}^{2,0}$}] {};

\node[new] at (intersection cs:
  first line={(c1) -- (c2)},
  second line={(g3) -- (g1)}) [label={[newlbl,label distance=-1ex]above left:$\vect{w}^{3,0}$}] {};

\node[new] at (intersection cs:
  first line={(c1) -- (c2)},
  second line={(g3) -- (g4)}) [label={[newlbl,label distance=-0.25ex]above right:$\vect{w}^{1,0}$}] {};
\end{tikzpicture}
\end{center}
\caption{The elementary step of the double description method}\label{fig:ddm}
\end{minipage}
\hfill
\begin{minipage}[t]{0.495\textwidth}
\vspace{0ex}
\begin{center}
\begin{tikzpicture}
[scale=0.8,>=triangle 45,
green/.style={circle,text=green!60!black,fill=green!60!black,inner sep=0pt,minimum size=1.25ex},
vtxlbl/.style={label distance=0pt,green!60!black},
red/.style={circle,text=red!80!black,fill=orange!30!red!90!black,inner sep=0pt,minimum size=1.25ex},
raylbl/.style={label distance=0pt,orange!30!red!90!black},
new/.style={circle,text=black,fill=blue!70,inner sep=0pt,minimum size=1.25ex},
newlbl/.style={label distance=0pt,blue!70}]
\equilateral{7}{90};

\barycenter{g}{\expo{0}}{\expo{0}}{\expo{0}};
\barycenter{gh}{\expo{1}}{\expo{0}}{\expo{0}};
\barycenter{h}{\expo{2}}{\expo{1}}{\expo{0}};

\barycenter{c1}{\expo{0}}{0}{0};
\barycenter{c2}{0}{\expo{0}}{\expo{0}};
\barycenter{c3}{0}{0}{\expo{0}};

\filldraw[blue!20,draw=none,fill opacity=0.3] (c1) -- (c2) -- (c3) -- cycle;
\draw[blue!40, ultra thick] (c1) -- (c2);

\draw[draw=black,ultra thick] (g) -- (gh) -- (h);
\node[red] at (h) [label={[raylbl]below:$\vect{w}$}] {};
\node[green] at (g) [label={[vtxlbl,label distance=-0.5ex]above left:$\vect{v}$}] {};
\node[new] at (gh) [label={[newlbl]left:$\vect{v}'$}] {};
\end{tikzpicture}
\end{center}
\caption{Illustration of Remark~\ref{rem:ddm_combination}}\label{fig:ddm_combination}
\end{minipage}
\end{figure}

\begin{example}\label{ex:ddm}
Figure~\ref{fig:ddm} provides an illustration of the elementary step on the cone defined in Example~\ref{ex:running} and the half-space given by the constraint $\vect{x}_2 \leq \vect{x}_3 + 5/2$ (depicted in light blue in Figure~\ref{fig:ddm}, while the set of elements activating the inequality is in darker blue). The three elements $\vect{v}^1$, $\vect{v}^2$, and $\vect{v}^3$ satisfy the constraint, while $\vect{v}^0$ does not. Their combinations are the elements $\vect{w}^{1,0}$, $\vect{w}^{2,0}$, and $\vect{w}^{3,0}$ respectively.
\end{example}

\begin{remark}\label{rem:ddm_combination}
The inductive scheme of the tropical double description method looks very similar to its classical counterpart, but they are distinguished by a minor difference: the combinations $(A_k \vect{w}) \vect{v} \mpplus (B_k \vect{v})\vect{w}$ with $A_k \vect{v} = B_k \vect{v}$ and $A_k \vect{w} > B_k \vect{w}$ do not appear in the classical case (see~\cite[Lemma~3]{FukudaProdon96}), while they are essential in the tropical setting. 

For instance, consider the cone of $\maxplus^3$ generated by the set $V$ consisting of the elements where $\vect{v} = (0,0,0)$ and $\vect{w} = (2,1,0)$ (in bold black in Figure~\ref{fig:ddm_combination}). Its intersection with the half-space $\{ (\vect{x}_1, \vect{x}_2, \vect{x}_3) \mid \vect{x}_2 \leq \vect{x}_3 \}$ (in light blue) is generated by a minimal set containing:
the vectors $\vect{v}$ (which activates the constraint $\vect{x}_2 \leq \vect{x}_3$) and $\vect{v}' = (2,1,1)$ (obtained by combining $\vect{v}$ and $\vect{w}$). Thus, the element $\vect{v}'$ cannot be dispensed with.
\end{remark}

\subsection{The tropical double description method algorithm}\label{subsec:resulting_algo}

Theorem~\ref{th:ddm_basis} and subsequently Theorem~\ref{th:ddm} may return non-extreme elements (see Example~\ref{ex:ddm} in which $\vect{w}^{2,0}$ and $\vect{w}^{3,0}$ are not extreme). If these redundant elements are not eliminated, the cardinality of the sets $V_k$ grows quadratically in the worst case at each step (because of the pairwise combinations in Theorem~\ref{th:ddm} of the $\vect{v}$ and $\vect{w}$). Hence the complexity of the inductive technique previously discussed is double exponential ($O(d^{2^p})$), both in time and space, which is clearly untractable. We propose to eliminate non-extreme elements at each step of the induction using the criterion based on directed hypergraphs, and the associated almost linear algorithm \Call{MaxScc}{}.

\begin{figure}
\begin{small}
\begin{algorithmic}[1]
\Procedure {ComputeExtRays}{$A,B,p$} \Lcomment{$A, B \in \maxplus^{p \times d}$}
  \If{$p = 0$} \Lcomment{Base case} \label{ce:base}
    \State \Return $\{\vect{\epsilon}^i\}_{i\in[d]}$
  \Else{} \Lcomment{Inductive case} \label{ce:inductive}
    \State $C \gets \left(\begin{smallmatrix}A_1 \\[-1ex] \vdots \\A_{p-1}\end{smallmatrix}\right)$, $D \gets \left(\begin{smallmatrix}B_1 \\[-1ex] \vdots \\B_{p-1}\end{smallmatrix}\right)$, $\vect{a} \gets A_p$, $\vect{b} \gets B_p$
    \State $V \gets \Call{ComputeExtRays}{C, D, p-1}$\label{ce:call}
    \State $V^\leq \gets \{ \vect{v}^i \in V \mid \vect{a} \vect{v}^i \leq \vect{b} \vect{v}^i  \}$, $V^> \gets \{ \vect{v}^j \in V \mid \vect{a} \vect{v}^j > \vect{b} \vect{v}^j  \}$, $W \gets V^\leq$ \label{ce:begin_elementary}
    \ForAll{$\vect{v}^i \in V^\leq$ and $\vect{v}^j \in V^>$}\label{ce:begin_loop}
      \State $\vect{w} \gets (\vect{a}\vect{v}^j)\vect{v}^i \mpplus (\vect{b} \vect{v}^i) \vect{v}^j$ 
      \State $\GG \gets \Call{BuildHypergraph}{\vect{w}, A, B}$\label{ce:build_hypergraph}
      \If{$\card{\Call{MaxScc}{\GG}} = 1$} \Lcomment{Extremality test}\label{ce:test}
	\State append $\mpnorm{\vect{w}}^{-1} \vect{w}$ to $W$ 
      \EndIf
    \EndFor\label{ce:end_loop} \label{ce:end_elementary}
  \EndIf
  \State \Return $W$
\EndProcedure
\end{algorithmic}
\end{small}
\caption{Implementation of the tropical double description method} \label{fig:computeextreme}
\end{figure}

The resulting algorithm \Call{ComputeExtRays}{} (Figure~\ref{fig:computeextreme}) provides the scaled representatives of the extreme rays of the cone $\CC$. The argument $p$ corresponds to the number of constraints of the system $A \vect{x} \leq B \vect{x}$. When $p = 0$, the cone coincides with $\maxplus^d$, and it is generated by the set  $\{\vect{\epsilon}^i\}_{i\in[d]}$. When $p > 0$, the system is split into the system $C \vect{x} \leq D \vect{x}$ formed by the first $(p-1)$ inequalities, and the last inequality $\vect{a} \vect{x} \leq \vect{b} \vect{x}$. The elements provided by Theorem~\ref{th:ddm} are computed from the set $V$ of extreme elements of the intermediate cone $\DD = \{\vect{x} \in \maxplus^d \mid C \vect{x} \leq D \vect{x} \}$. The set $W$ is used to store the extreme rays of $\CC$. The extremality test is evaluated at Lines~\lineref{ce:build_hypergraph}-\lineref{ce:test}. First, the tangent hypergraph $\GG(\vect{w}, \CC)$ is computed thanks to a function \Call{BuildHypergraph}{}. Then, the function \Call{MaxScc}{} returning the set of the maximal \scc{}s of the hypergraph is called. 
If the test succeeds, the element $\vect{w}$ is first normalized into the scaled element $\mpnorm{\vect{w}}^{-1} \vect{w}$, and then appended to the set $W$. 

Observe that the extremality test is applied only to the elements associated with the combinations $(\vect{a}\vect{v}^j)\vect{v}^i \mpplus (\vect{b} \vect{v}^i) \vect{v}^j$, and not to the elements $\vect{v} \in V^\leq$ which satisfy $\vect{a} \vect{v} \leq \vect{b} \vect{v}$. Indeed, each element $\vect{v} \in V^\leq$ is extreme in the cone $\DD$, and subsequently in the cone $\CC$, since $\CC \subset \DD$.

\subsubsection*{Complexity analysis.}

Each operation in $\maxplus$ is supposed to take a unit time. We use hash sets to encode subsets of $\maxplus^d$. A \emph{hash set} is a hash table which maps keys to a fixed value (chosen arbitrarily, for instance $\mathit{Nil}$). The keys stored in the hash table correspond to the elements of the represented set. The amortized time complexity of adding, searching, and removing an element in the set is bounded by the complexity of hashing a vector of $\maxplus^d$, which is supposed to be $O(d)$. 

We first study the complexity of the \emph{inductive step}, which refers to the set of operations performed since the last call to \Call{ComputeExtRays}{} (Lines~\lineref{ce:begin_elementary} to~\lineref{ce:end_elementary}). Starting from the last intermediate generating set $V$, it consists in
\begin{inparaenum}[(i)]
\item computing the set given in~\eqref{eq:ddm_basis}, and 
\item eliminating non-extreme combinations.
\end{inparaenum}
Its complexity can be precisely characterized in terms of the size of $V$. It can be verified that it is dominated by the complexity of the extremality tests performed in the loop from Lines~\lineref{ce:begin_loop} to~\lineref{ce:end_loop}. Each test requires to build a hypergraph $\GG$ (Line~\lineref{ce:build_hypergraph}). This operation can be done in linear time in its size, which is in $O(p d)$. According to Theorem~\ref{th:complexity}, $\Call{MaxScc}{\GG}$ is executed in time $O(\size(\GG) \alpha(d)) = O(p d \alpha(d))$. The loop is iterated $O(\card{V}^2)$ in the worst case, so that the following statement holds: 
\begin{proposition}\label{prop:ce_inductive_complexity}
The worst case time complexity of the inductive step in \Call{ComputeExtRays}{} is $O(p d \alpha(d) \card{V}^2)$.
\end{proposition}
We also stress that the inductive step is optimal in terms of space complexity, since a non-extreme element is never stored in the resulting set $W$, even temporarily. It follows that its space complexity is bounded by $O(d \max(\card{V},\card{W}))$.

\begin{remark}\label{remark:argmax}
Observe that the construction of the hypergraph $\GG$ (Line~\lineref{ce:build_hypergraph}) can be optimized by maintaining some extra information for each element of the intermediate set $V$.

Indeed, consider a tropical linear form $\vect{c} \in \maxplus^{1 \times d}$ and a non-null combination $\vect{v} = \lambda \vect{x} \mpplus \mu \vect{y}$ of two elements $\vect{x} , \vect{y} \in \maxplus^d$. The set $\argmax(\vect{c} \vect{v})$ can be computed efficiently from the sets $\argmax(\vect{c} \vect{x})$ and $\argmax(\vect{c} \vect{y})$:
\begin{equation}
\argmax(\vect{c} \vect{v}) = 
\begin{cases}
\argmax(\vect{c} \vect{x}) & \text{if }\lambda (\vect{c} \vect{x}) > \mu (\vect{c} \vect{y}), \\
\argmax(\vect{c} \vect{y}) &  \text{if }\lambda (\vect{c} \vect{x}) < \mu (\vect{c} \vect{y}), \\
\argmax(\vect{c} \vect{x}) \cup \argmax(\vect{c} \vect{y}) & \text{otherwise}.
\end{cases} \label{eq:argmax}
\end{equation}
The value of $\vect{c} \vect{v}$ can be computed in $O(1)$ time from $\vect{c} \vect{x}$ and $\vect{c} \vect{y}$ using the same idea.

Now, let $\vect{v}$ be an element returned by \Call{ComputeExtRays}{$A, B, p$}. Using~\eqref{eq:argmax}, the list of the tuples $((A_k \vect{v}, \argmax(A_k \vect{v})), (B_k \vect{v}, \argmax(B_k \vect{v})))$ ($k\in[p]$) can be propagated by induction during the execution of \Call{ComputeExtRays}{$A, B, p$}. In practice, we have observed that this optimization considerably speeds up the computation of the associated hypergraph.
\end{remark}

The overall complexity of the algorithm \Call{ComputeExtRays}{} depends on the maximal size of the sets $V_k$ ($0 \leq k \leq p-1$) returned in the intermediate steps:
\begin{proposition}\label{prop:ce_overall_complexity1}
The worst case time complexity of the \Call{ComputeExtRays}{} algorithm is bounded by 
\[
O(p^2 d \alpha(d) V^2_{\max}),
\]
where $V_{\max}$ is the maximal cardinality of the sets $V_k$ for $k = 0,\dots,p-1$.
\end{proposition}

In classical geometry, the upper bound theorem of McMullen~\cite{mcmullen70} shows that the maximal number of extreme points of a convex polytope in $\real^d$ defined by $p$ inequality constraints is equal to 
\[
U(p,d) \defi 
\begin{dcases}
\binom{p-\floor{d/2}}{\floor{d/2}} + \binom{p-\floor{d/2}-1}{\floor{d/2}-1} & \text{for $d$ even,} \\
2\binom{p-\floor{d/2}-1}{\floor{d/2}} &\text{for $d$ odd.}
\end{dcases}
\]
The polars of the \emph{cyclic polytopes} (see~\cite{ziegler98}) are known to reach this bound. In the tropical setting, a recent work of Allamigeon, Gaubert, and Katz~\cite{AllamigeonGaubertKatzJCTA2011} proves that the number of extreme rays of a tropical polyhedral cone $\CC$ in $\maxplus^d$ defined by $p$ inequalities is bounded by a similar quantity.
\begin{theorem}[{\cite{AllamigeonGaubertKatzJCTA2011}}]\label{th:upperbound}
The number of extreme rays of a tropical cone in $\maxplus^d$ defined as the intersection of $p$ tropical half-spaces cannot exceed $U(p+d,d-1)$.
\end{theorem}
The bound of Theorem~\ref{th:upperbound} is shown in~\cite{AllamigeonGaubertKatzJCTA2011} to be asymptotically tight for a fixed $p$, as $d$ tends to infinity, being approached by the \emph{signed cyclic tropical polyhedral cones}, which are a tropical generalization of the (polar of) the cyclic polytope, taking into account a sign pattern. The bound is believed not to be tight for a fixed $d$, as $p$ tends to infinity,  because the growth of the number of extreme rays for the model of signed cyclic polyhedral cones is too slow. Finding the optimal bound is an open problem. 

By combining Proposition~\ref{prop:ce_overall_complexity1} and Theorem~\ref{th:upperbound}, we readily get the following upper bound on the complexity of \Call{ComputeExtRays}{}:
\begin{corollary}\label{cor:ce_overall_complexity2}
The worst-case time complexity of the \Call{ComputeExtRays}{} algorithm is bounded by $O\bigl(p^2 d \alpha(d) \cdot (U(p+d-1,d-1))^2\bigr)$. 
\end{corollary}
The asymptotic behavior of the bound of Corollary~\ref{cor:ce_overall_complexity2} is the following one:
\begin{equation}
O\biggl(p^2 d \alpha(d) \cdot \Bigl(1 + \frac{p}{\lfloor (d-1)/2\rfloor}\Bigr)^{2\lfloor \frac{d-1}{2} \rfloor}  \Bigl(1 + \frac{\lfloor (d-1)/2\rfloor}{p}\Bigr)^{2p}\biggr) \quad \text{when } p+d \gg 1. \label{eq:complexity_bound}
\end{equation}
In particular, the bound~\eqref{eq:complexity_bound} is dominated by $O(p^2 d \alpha(d) \cdot (e(1+\frac{M}{m})^{2m}))$, where $m$ and $M$ are respectively the minimum and the maximum of $p$ and $\lfloor (d-1)/2\rfloor$. The algorithm \Call{ComputeExtRays}{} is therefore polynomial time as soon as one of the parameters $d$ or $p$ is constant. 
In general, we shall keep in mind that the quality of the bounds given above directly depends on the quality of the upper bound of Theorem~\ref{th:upperbound}. Since the latter
may not be tight for certain asymptotic regimes of the $(p,d)$ parameters,
the former may give a loose overestimate of the complexity
of the algorithm \Call{ComputeExtRays}{}.

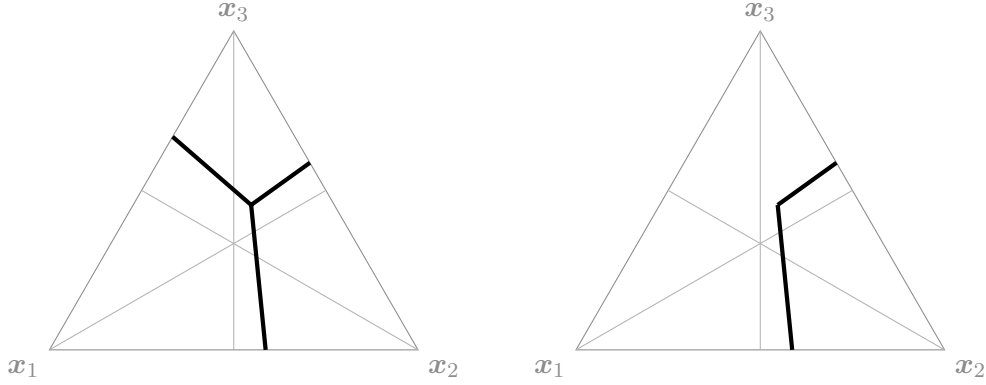
\begin{figure}
\begin{center}
\begin{tikzpicture}[scale=0.7,>=triangle 45]
\begin{scope}
\equilateral{7}{90};

\barycenter{c}{\expo{0}}{\expo{0.5}}{\expo{1}};
\barycenter{cxy}{\expo{0}}{\expo{0.5}}{0};
\barycenter{cyz}{0}{\expo{0.5}}{\expo{1}};
\barycenter{cxz}{\expo{0}}{0}{\expo{1}};

\draw[draw=black,ultra thick] (c) -- (cxy) (c) -- (cyz) (c) -- (cxz);
\end{scope}
\begin{scope}[shift={(10,0)}]
\equilateral{7}{90};

\barycenter{c}{\expo{0}}{\expo{0.5}}{\expo{1}};
\barycenter{cxy}{\expo{0}}{\expo{0.5}}{0};
\barycenter{cyz}{0}{\expo{0.5}}{\expo{1}};
\barycenter{cxz}{\expo{0}}{0}{\expo{1}};

\draw[draw=black,ultra thick] (c) -- (cxy) (c) -- (cyz);
\end{scope}
\end{tikzpicture}
\end{center}
\caption{The tropical hyperplane ``$\max(\vect{x}_1,1/2+\vect{x}_2,1+\vect{x}_3)$ attained at least twice'' (left), and the signed tropical hyperplane $1/2+\vect{x}_2 = \max(\vect{x}_1,1+\vect{x}_3)$ (right)}\label{fig:hyperplane}
\end{figure}

\subsection{Variants of the algorithm with other kinds of constraints}\label{subsec:hyperplanes}

Our algorithm defined in Section~\ref{subsec:resulting_algo} can be generalized to handle systems including tropical linear constraints other than inequalities.

\subsubsection{Tropical hyperplanes} 

Tropical geometry originates when looking at classical objects with logarithmic glasses or valuations. Let $k$ denotes the field of complex Puiseux series in an indeterminate $t$, equipped with the valuation $v$ which takes the opposite of the smallest exponent arising in a series. Then, a tropical linear space may be defined as the closure of the image of a linear space over $k$ by the map which applies the valuation $v$ to each coordinate, see~\cite{SpeyerSturmfels04,RGST} for more information. In particular,
consider the hyperplane in $k^d$ defined by the equation
$\sum_{i = 1}^d a_i\vect{x}_i = 0$, with $a_i\in k$.
Then, a theorem of Kapranov characterizing more generally the nonarchimedian amoebas of hypersurfaces (see~\cite{kapranov}) shows that the closure of the image of this hyperplane by the valuation coincides with the set of points $\vect{y} \in \maxplus^d$ such that the maximum in the expression $\max_{1 \leq i \leq d} v(a_i)+ \vect{y}_i$ is attained at least twice. Such a set is known as a \emph{tropical hyperplane}~\cite{RGST,DS}. Tropical hyperplanes form a subclass of tropical polyhedral cones of $\maxplus^d$, see Figure~\ref{fig:hyperplane}. 

Tropical hyperplanes can be handled in the elementary step of the double description method thanks to the following result:
\begin{theorem}\label{th:ddm_basis_hyperplane}
Let $\CC \subset \maxplus^d$ be a closed tropical cone generated by a set $V$ of elements of $\maxplus^d$, and let $\HH = \{ \vect{x} \in \maxplus^d \mid \vect{c} \vect{x} = \max_{i \in [d]} (\vect{c}_i + \vect{x}_i) \text{ is attained at least twice} \}$ be a tropical hyperplane ($\vect{c} \in \maxplus^{1 \times d}$). Then the cone $\CC \cap \HH$ is generated by:
\begin{equation}
\begin{aligned}
& \phantom{{}\cup{}} \{ \vect{v} \in V \mid \vect{c} \vect{v} \text{ is attained at least twice} \} \\
& {}\cup \{ (\vect{c} \vect{w})\vect{v} \mpplus (\vect{c} \vect{v}) \vect{w} \mid \vect{v}, \vect{w} \in V \text{, } \vect{c} \vect{v} \text{ attained only once and } \vect{c} \vect{w} \text{ at least twice} \} \\
& {}\cup \{ (\vect{c} \vect{w})\vect{v} \mpplus (\vect{c} \vect{v}) \vect{w} \mid \vect{v}, \vect{w} \in V \text{, } \vect{c} \vect{v} \text{ and } \vect{c} \vect{w} \text{ are attained only once and at distinct indices} \}.
\end{aligned} \label{eq:ddm_basis_hyperplane}
\end{equation}
\end{theorem}

\begin{proof}
If $V'$ is the set of vectors given in~\eqref{eq:ddm_basis_hyperplane}, then we clearly have $V' \subset \CC \cap \HH$ so that $\mpcone(V') \subset \CC \cap \HH$. 

Conversely, supposing $\vect{x} \in \CC \cap \HH$, then by Minkowski-Carath{\'e}odory theorem on $\CC = \mpcone(V)$, we have $\vect{x} = \mpsum_{i = 1}^d \lambda_i \vect{v}^i$ with $\vect{v}^i \in V$ ($i = 1,\dots,d$). Let $I, J \subset \oneto{d}$ such that for all $i \in I$, $\vect{c} \vect{v}^i$ is attained at least twice, and for every $j \in J$, $\vect{c} \vect{v}^j$ is attained only once ($I \cap J = \emptyset$). 

We know that the maximum $\mpsum_{i \in I} \lambda_i (\vect{c} \vect{v}^i) \mpplus \mpsum_{j \in J} \lambda_j (\vect{c} \vect{v}^j)$ is reached at least twice. Let $\kappa = \vect{c} \vect{x}$. If $\kappa = \mpzero$, then necessarily for any $j \in J$, $\lambda_j (\vect{c} \vect{v}^j) = \mpzero$ hence $\lambda_j = \mpzero$ (if not, $\vect{c} \vect{v}^j = \mpzero$ would be attained more than once). Thus $\vect{x} \in \mpcone(V')$ obviously. 

Now, suppose that $\kappa > \mpzero$. We distinguish two cases:
\begin{enumerate}[(i)]
\item suppose that there exists $i_0 \in I$ such that $\lambda_{i_0} (\vect{c} \vect{v}^{i_0}) = \kappa$. Then:
\begin{align*}
\vect{x} & = \mpsum_{i \in I} \lambda_i \vect{v}^i \mpplus \mpsum_{j \in J} \lambda_j \vect{v}^j \\
& = 
\mpsum_{i \in I} \lambda_i \vect{v}^i \mpplus \kappa^{-1} \biggl(\mpsum_{j \in J} \lambda_j (\vect{c} \vect{v}^j)\biggr) \lambda_{i_0} \vect{v}^{i_0}  \mpplus \kappa^{-1} \mpsum_{j \in J} \bigl(\lambda_{i_0} (\vect{c} \vect{v}^{i_0})\bigr) \lambda_j \vect{v}^j 
&& \text{as } \kappa \geq \mpsum_{j \in J} \lambda_j (\vect{c} \vect{v}^j)
\\
& = \mpsum_{i \in I} \lambda_i \vect{v}^i \mpplus \kappa^{-1} \mpsum_{j \in J} \lambda_{i_0} \lambda_j \bigl((\vect{c} \vect{v}^j) \vect{v}^{i_0} \mpplus  (\vect{c} \vect{v}^{i_0}) \vect{v}^j\bigr) .
\end{align*}

\item otherwise, for all $i \in I$, $\lambda_i (\vect{c} \vect{v}^i) < \vect{c} \vect{x}$. In this case, the maximum $\vect{c} \vect{x}$ is necessarily attained by two distinct terms $\lambda_{j_1} (\vect{c} \vect{v}^{j_1})$ and $\lambda_{j_2} (\vect{c} \vect{v}^{j_2})$, with $j_1, j_2 \in J$, $j_1 \neq j_2$, and if $k_1$ and $k_2$ are respectively the arguments of the maxima $\vect{c} \vect{v}^{j_1}$ and $\vect{c} \vect{v}^{j_2}$, we have $k_1 \neq k_2$. Let $J_l = \{ j \in J \mid k_l \not \in \argmax(\vect{c} \vect{v}^j)\}$ for $l = 1, 2$. Note that $J_1 \cup J_2 = J$. Then:
\begin{align*}
\vect{x} & = \mpsum_{i \in I} \lambda_i \vect{v}^i  \mpplus \mpsum_{j \in J} \lambda_j \vect{v}^j \\ 
& = \mpsum_{i \in I} \lambda_i \vect{v}^i  \mpplus \kappa^{-1} \biggl(\mpsum_{j \in J_1} \lambda_j (\vect{c} \vect{v}^j) \biggr) \lambda_{j_1} \vect{v}^{j_1}  \mpplus  \kappa^{-1} \mpsum_{j \in J_1} \bigl( \lambda_{j_1} (\vect{c} \vect{v}^{j_1}) \bigr) \lambda_j\vect{v}^j \mpplus  \kappa^{-1} \biggl(\mpsum_{j \in J_2} 
\lambda_j (\vect{c} \vect{v}^j) \biggr) \lambda_{j_2} \vect{v}^{j_2} \\
& \phantom{{}={}}
{}\mpplus  \kappa^{-1} \mpsum_{j \in J_2} \bigl( \lambda_{j_2} (\vect{c} \vect{v}^{j_2}) \bigr) \lambda_j\vect{v}^j \ \ \ \ \ \ \text{since } \kappa = \lambda_{j_1} (\vect{c} \vect{v}^{j_1}) = \lambda_{j_2} (\vect{c} \vect{v}^{j_2}), \text{ and } \kappa \geq \mpsum_{j \in J} \lambda_j (\vect{c} \vect{v}^j)
\\
& = \mpsum_{i \in I} \lambda_i \vect{v}^i   \mpplus  \kappa^{-1} \mpsum_{j \in J_1} 
\lambda_j \lambda_{j_1} \bigl( (\vect{c} \vect{v}^j) \vect{v}^{j_1} \mpplus (\vect{c} \vect{v}^{j_1}) \vect{v}^j \bigr) \mpplus  \kappa^{-1} \mpsum_{j \in J_2} 
\lambda_j \lambda_{j_2} \bigl( (\vect{c} \vect{v}^j) \vect{v}^{j_2} \mpplus (\vect{c} \vect{v}^{j_2}) \vect{v}^j \bigr).
\end{align*}
\end{enumerate}
In both cases, $\vect{x} \in \mpcone(V')$, which completes the proof.
\end{proof}

The extremality criterion of Theorem~\ref{th:main} can be extended to systems containing tropical hyperplane constraints. Every hyperplane $\{ \vect{x} \in \maxplus^d \mid \vect{c} \vect{x} \text{ is attained at least twice} \}$ generates the hyperarcs $(\argmax(\vect{c} \vect{v}) \setminus \{ i \}, \{ i \})$ (for each $i \in \argmax(\vect{c}\vect{v})$) in the directed hypergraph $\GG(\vect{v},\CC)$. This results from the fact that the hyperplane can be equivalently expressed as the set of the solutions of the following system:
\begin{equation}
\vect{c}_i \vect{x}_i \leq \mpsum_{\substack{1 \leq j \leq d \\ j \neq i}} \vect{c}_j \vect{x}_j \qquad \text{for } i = 1,\dots,d \label{eq:hyperplane_system}
\end{equation}
and that~\eqref{eq:hyperplane_system} is active on $\vect{v}$ if, and only if, $i \in \argmax(\vect{c} \vect{v})$.

Thus, the extreme rays of the intersection of $p$ tropical hyperplanes can be determined in time $O((p d)^2 \alpha(d) (V'_{\max})^2)$, where $V'_{\max}$ is the maximal size of the sets arising in the intermediate steps of the induction. In comparison, by expanding hyperplanes to a collection of $p d$ half-spaces using~\eqref{eq:hyperplane_system}, the complexity of the algorithm~\Call{ComputeExtRays}{} is $O(p^2 d^3 \alpha(d) V_{\max})$. Since $V'_{\max} \leq V_{\max}$, the variant presented here improves time bounds by a factor $d V_{\max}/V'_{\max}$.

\subsubsection{Signed tropical hyperplanes} 

Another noticeable case of tropical polyhedral cones consists of \emph{signed tropical hyperplanes}, which are sets of points satisfying an equality $\vect{a} \vect{x} = \vect{b} \vect{x}$, where the support of the vectors $\vect{a}$ and $\vect{b}$ are disjoint (see~\cite{Plus,AGG08b}). They correspond to subsets of (non-signed) hyperplanes (Figure~\ref{fig:hyperplane}). The elementary step of the double description method can be extended to such sets, as follows:
\begin{theorem}\label{th:ddm_basis_eq}
Let $\CC \subset \maxplus^d$ be a closed tropical cone generated by a set $V$ of elements of $\maxplus^d$, and let $\vect{a}, \vect{b} \in \maxplus^{1 \times d}$. Then the cone $\CC \cap \{ \vect{x} \in \maxplus^d \mid \vect{a} \vect{x} = \vect{b} \vect{x} \}$ is generated by the following set:
\begin{align*}
& \phantom{{}\cup{}} \{ \vect{v} \in V \mid\vect{a} \vect{v} = \vect{b} \vect{v} \} \cup \{ (\vect{a} \vect{w}) \vect{v} \mpplus (\vect{b} \vect{v}) \vect{w} \mid \vect{v},\vect{w} \in V \text{, } \vect{a} \vect{v} < \vect{b} \vect{v} \text{ and } \vect{a} \vect{w} > \vect{b} \vect{w} \} 
 \\
& {}\cup \{ (\vect{a} \vect{w}) \vect{v} \mpplus (\vect{b} \vect{v}) \vect{w}) \mid \vect{v},\vect{w} \in V \text{, } \vect{a} \vect{v} < \vect{b} \vect{v} \text{ and } \vect{a} \vect{w} = \vect{b} \vect{w} \} \\
& {}\cup \{ (\vect{a} \vect{w}) \vect{v} \mpplus (\vect{b} \vect{v}) \vect{w} \mid \vect{v},\vect{w} \in V \text{, } \vect{a} \vect{v} = \vect{b} \vect{v} \text{ and } \vect{a} \vect{w} > \vect{b} \vect{w} \}.
\end{align*}
\end{theorem}

\begin{proof}
Straightforward from two successive applications of Theorem~\ref{th:ddm_basis} on the inequalites $\vect{a} \vect{x} \leq \vect{b} \vect{x}$ and $\vect{a} \vect{x} \geq \vect{b} \vect{x}$.
\end{proof}

The extremality criterion of Theorem~\ref{th:main} can also be generalized to signed hyperplanes, by introducing two symmetric hyperarcs per equality in the tangent directed hypergraph: for an equality $\vect{a} \vect{x} = \vect{b} \vect{x}$, the tangent hypergraph $\GG(\vect{v},\CC)$ at the element $\vect{v}$ will contain the hyperarcs $(\argmax(\vect{a} \vect{v}), \argmax(\vect{b} \vect{v}))$ and $(\argmax(\vect{b} \vect{v}), \argmax(\vect{a} \vect{v}))$. 

\section{Arrangements of tropical hyperplanes}\label{sec:arrangement}

In this section, for the sake of comparison, we present an alternative approach to the problem of computing the extreme rays of a tropical cone described as the intersection of half-spaces in general position. This approach relies on arrangements of signed tropical hyperplanes. 

We consider the case in which the tropical cone $\CC$ is defined as the set of the solutions $\vect{x} \in \maxplus^d$ of the system of inequalities $A \vect{x} \leq B \vect{x}$ ($A, B \in \maxplus^{p \times d}$). We suppose, without loss of generality, that for all $k \in \oneto{p}$, the supports of the $k$-th rows of $A$ and $B$ are disjoint. For each half-space $\{ \vect{x} \in \maxplus^d \mid A_k \vect{x} \leq B_k \vect{x} \}$, we introduce the \emph{associated signed hyperplane}, denoted by $H_k$ and defined as the set of the elements $\vect{x} \in \maxplus^d$ satisfying the equality $A_k \vect{x} = B_k \vect{x}$. 

Following~\cite{DSS,RGST}, a matrix $M \in \maxplus^{k \times k}$ is said to be \emph{tropically non-singular} if, and only if, the tropical permanent
\begin{equation}
\tper M = \mpplus_{\sigma \in S_k} M_{1\sigma(1)}\dots M_{k\sigma(k)} \label{eq:tper}
\end{equation}
is not null, and the maximum in~\eqref{eq:tper} is reached at precisely one permutation $\sigma$ in the symmetric group $S_k$. 
In this section, we will assume that the half-spaces defining the cone $\CC$ are in \emph{general position}, meaning that every square submatrix of $A \mpplus B$ is (tropically) non-singular.

A finite set of signed hyperplanes constitutes an \emph{arrangement}.
In this setting, a \emph{vertex} of the arrangement will refer to the intersection of some of the hyperplanes when this intersection is reduced to a ray $\maxplus \vect{x}$ (\ie{}, a point, written in homogeneous coordinates).
\begin{proposition}\label{prop:arrangement}
When the half-spaces defining the cone $\CC$ are in general position, the extreme rays of $\CC$ are vertices of the arrangement formed  by the signed hyperplanes $H_k$ ($1 \leq k \leq p$) and the hyperplanes $Z_j = \{\vect{x} \in \maxplus^d \mid \vect{x}_j = \mpzero \}$ ($1 \leq j \leq d$). 
\end{proposition}

Before proving Proposition~\ref{prop:arrangement}, we first recall that the max-plus Cramer theorem~\cite{Plus,RGST,AGG08b} shows that the vertices of the arrangement of Proposition~\ref{prop:arrangement} are precisely given as the 
non-trivial
intersections of $(d-1)$ signed hyperplanes. It also provides a constructive method to determine the vertices by using Cramer permanents. 
\begin{proposition}[{Corollary of~\cite{Plus},\cite[\S~5]{RGST},\cite[Th.~6.6]{AGG08b}}]\label{prop:cramer}
Given $n \in \oneto{d}$, let $A'$ (resp.\ $B'$) the sub-matrix formed by the first $(n-1)$ rows and $n$ columns of $A$ (resp.\ $B$). Let $C_i$ be the matrix of size $(n-1) \times (n-1)$ obtained from $A' \mpplus B'$ by deleting the $i$-th column ($i \in \oneto{n}$). Let $\vect{x} \in \maxplus^d$ be the vector of support $\oneto{n}$ defined by $\vect{x}_i = \tper C_i$ for all $i \in \oneto{n}$. 

Then the intersection of the hyperplanes $H_k$ ($1 \leq k \leq n-1$) and $Z_j$ ($n+1 \leq j \leq d$) is either empty, or reduced to the ray $\maxplus \vect{x}$.
\end{proposition}
These $d$ Cramer permanents $\tper C_k$ can be naively computed by solving $d$ assignment problems, leading to a time complexity $O(d^4)$. However, as remarked in~\cite{RGST}, all the permanents can be determined (up to a multiplicative constant) as the optimal solution of a single transportation problem. This allows to determine the vertex of the arrangement in time $O(d^3)$. Alternatively, 
all the Cramer permanents can be determined by solving a single optimal assignment problem, and
then by applying a variant of the Jacobi algorithm in~\cite{Plus}, which also
gives a $O(d^3)$ algorithm.
Note also, although that it will not be needed here, that the emptyness of the intersection can be checked a priori by inspecting the parity of the optimal permutations in the tropical Cramer permanents~\cite{Plus}, \cite[Th.~6.4]{AGG08b}. 
Finally, observe that Proposition~\ref{prop:cramer} applies to the intersection of hyperplanes $H_k$ ($k \in K$), and $Z_j$ ($j \in J$), with $\card{J} + \card{K} = d-1$, up to permuting hyperplanes and coordinates.

\begin{proof}[\proofname{} of Proposition~\ref{prop:arrangement}]
Given a representative $\vect{v}$ of an extreme ray of $\CC$, we are going to show that at least $(n-1)$ inequalities among the system $A \vect{x} \leq B \vect{x}$ are active at the point $\vect{v}$, where $n$ is the cardinality of the support of $\vect{v}$. 

Without loss of generality, we will assume that $\supp(\vect{v}) = \oneto{n}$, with $n \leq d$. Let $\vect{w}$ be the vector of $\maxplus^n$ reduced to the first $n$ coordinates of $\vect{v}$. Applying Proposition~\ref{prop:support_extremality}, the element $\vect{w}$ is extreme in the cone $\CC'$ defined the system $A' \vect{y} \leq B' \vect{y}$ over $\vect{y} \in \maxplus^n$, where $A'$ and $B'$ respectively correspond to the first $n$ columns of $A$ and $B$. By Theorem~\ref{th:main}, let $i$ be a node of the greatest \scc{} of the tangent directed hypergraph $\GG(\vect{w},\CC')$. Any other node $j \in \oneto{n} \setminus \{i\}$ of $\HH(\vect{w},\CC')$ has to reach $i$. Hence, there exists a hyperarc $a_{k_j}$ such that $T(a_{k_j})$ is reduced to the singleton $\{j\}$ (if not, $j$ would not reach any node except itself). Each hyperarc corresponds to an inequality active at $\vect{w}$, and subsequently at $\vect{v}$.

The extreme ray $\maxplus \vect{v}$ is therefore included into the intersection of the signed hyperplanes $H_{k_j}$ with $j \in \oneto{n} \setminus \{i\}$ and $Z_{n+1}, \dots, Z_d$. By Proposition~\ref{prop:cramer} and the subsequent discussion, we deduce that the intersection of the hyperplanes is reduced to the ray $\maxplus \vect{v}$. The latter is consequently a vertex of the arrangement.
\end{proof}

\begin{remark}\label{remark:vertex_not_extreme}
Observe that the converse of Proposition~\ref{prop:arrangement} does not hold, in the sense that not every vertex of the arrangement belonging to the cone $\CC$ is an extreme ray. For instance, for $d = 3$, consider the cone $\CC$ defined as the intersection of the half-spaces associated with the inequalities $\vect{x}_1 \mpplus (-1) \vect{x}_2 \leq \vect{x}_3$ and $(-1)\vect{x}_1 \mpplus \vect{x}_2 \leq \vect{x}_3$. These two half-spaces are in general position, and the intersection of the associated hyperplanes are reduced to the ray $\maxplus (0,0,0)$, which is not extreme in the cone $\CC$.
\end{remark}

\begin{figure}
\begin{algorithmic}[1]
\State $V \gets \emptyset$
\ForAll{subset $S$ of $(d-1)$ hyperplanes among the $Z_j$ ($j \in \oneto{d}$) and $H_k$ ($k \in \oneto{p}$)}
  \If{the intersection $\cap_{H \in S} H$ is a ray $\maxplus \vect{v}$ such that $\vect{v} \in \CC$}
  \State $\GG \gets \Call{BuildHypergraph}{\vect{v}, A, B}$
  	\If{$\card{\Call{MaxScc}{\GG}} = 1$} 
	\State append $\mpnorm{\vect{v}}^{-1} \vect{v}$ to $V$
	\EndIf
  \EndIf
\EndFor
\State \Return $V$
\end{algorithmic}
\caption{Computing the extreme rays from the vertices of signed hyperplane arrangement}\label{fig:arrangement_algo}
\end{figure}

Proposition~\ref{prop:arrangement} naturally leads to the idea of determining the extreme rays of the cone $\CC$ by enumerating of all the vertices of the arrangement, then keeping only the vertices belonging to the cone, and eliminating the non-extreme ones using the characterization of Theorem~\ref{th:main}. This algorithm is provided in Figure~\ref{fig:arrangement_algo}. Its complexity is given by the following result:
\begin{theorem}
When the half-spaces defining the tropical cone $\CC$ are in general position, the extreme rays of $\CC$ can be determined in time $O\bigl((p d \alpha(d) + d^3) \cdot \binom{p+d}{d-1}\bigr)$.
\end{theorem}
The complexity of the algorithm is thus asymptotically given by:
\[
O\biggl((p d \alpha(d) + d^3) \cdot \Bigl(1 + \frac{p+1}{d-1}\Bigr)^{d-1} \Bigl(1 + \frac{d-1}{p+1}\Bigr)^{p+1}\biggr) \quad \text{when } p+d \gg 1.
\]
In theory, this improves the worst-case bound of the algorithm~\Call{ComputeExtRays}{} given in Corollary~\ref{cor:ce_overall_complexity2} in some cases. Nevertheless, in practice, the algorithm of Figure~\ref{fig:arrangement_algo} appears to be of little use, since the worst case execution time is essentially always achieved. The leading term $\binom{p+d}{d-1}$ of the complexity bound is indeed a lower bound on the execution time of the algorithm (every subset $S$ of $(d-1)$ hyperplanes is examined). In contrast, the algorithm \Call{ComputeExtRays}{} takes advantage of the fact that $V_{\max}$ is in general much smaller than the upper bound of Theorem~\ref{th:upperbound}. 

Furthermore, the algorithm presented here can only be applied to the general position setting, contrary to \Call{ComputeExtRays}{}. 
Given an arbitrary polyhedral cone $\CC$, this is certainly possible to execute the algorithm of Figure~\ref{fig:arrangement_algo} on an intersection of half-spaces in general position approximating $\CC$. However, this approach may provide many more rays than the cone $\CC$ actually has. Indeed, the approximation of a polyhedral cone $\CC$ by a sequence of polyhedral cones $\CC(m)$ defined by half-spaces in general position is discussed in Section~5 of \cite{AllamigeonGaubertKatzJCTA2011}. The proof of Theorem~8 there deals with the case where $\CC\subset \CC(m)$ holds for all $m$. It shows that any accumulation point of a family of representatives of extreme rays of $\CC(m)$, as $m\to\infty$, is a generating family of $\CC$. However, in general the generating family obtained in this way contains many redundant generators.

\section{Comparison with alternative approaches and experimental results}\label{sec:comparison}\label{subsec:existing_approaches}

\subsection{Existing incremental algorithms in the tropical setting.}\label{subsec:existing_incremental}
Butkovi\v{c} and Hegedus~\cite{butkovicH} 
were the first to establish the existence of a finite generating family
for the set of solutions of $A\vect{x}\leq B\vect{x}$ (see also~\cite{ButkovicBook} for a recent account), but their method was not
intended to be algorithmically efficient. Unlike the present approach,
their elementary step leads to an algorithm which is not tail-recursive, 
and in which it is not natural to incorporate an incremental elimination
of redundant generators. This provides a method with a double exponential complexity.

The principle of the approach implemented by the second author in the Maxplus toolbox of {\sc Scilab}~\cite{toolbox}, and refined in our previous work~\cite{AGG08}, is similar to the one of \Call{ComputeExtRays}{}. However, it uses a much less efficient elimination of non-extreme vectors, which does not take the external representation of the cone into account. This elimination relies on a characterization in terms of set covers due to Vorobyev and Zimmermann, or equivalently on residuation (see~\cite{vorobyev67,Zimmermann76} and~\cite{agk04} for a recent overview). Its time complexity is $O(d \card{V}^2)$, where $V$ is the set of extreme rays of the previous intermediate cone $\DD$. Subsequently, the total complexity of the inductive step is $O(d \card{V}^4)$.

Butkovi\v{c}, Schneider, and Sergeev~\cite{BSS} proposed
a characterization of extreme points in terms of minimal
elements of a given type (see Remark~\ref{remark:minimality} above) in order to compute a minimal generating family of a polyhedral cone given
by a set of rays. This characterization reduces the latter problem
to computing the set of (Pareto) minimal vectors of
a given set of $k$ vectors in dimension $d$, which can be solved in time time $O(k(\log_2 k)^{d-3})$ (for $d \geq 4$) using~\cite{kunglucciopreparata}. Using this approach as a replacement of Theorem~\ref{th:main}, this leads to a variant in which the inductive step has a complexity of $O(d\card{V}^22^{d-3}(\log_2\card{V})^{d-3})$ (the algorithm of ~\cite{kunglucciopreparata} must be called $d$ times, on a set of $O(\card{V}^2)$ vectors).

\subsubsection*{Time complexities}

The complexities of the different approaches are compared in Table~\ref{tab:comparison}. Recall that $V_{\max}$ denotes the maximal cardinality of the intermediate generating sets. As shown in~\cite{AllamigeonGaubertKatzJCTA2011}, there exists instances
in which $V_{\max}$ is of order $(p+d)^{(d-1)/2}$, for $d\gg p$,
other instances in which it is of order $p2^{d-2}$ for $p\gg d$,
and in general, experiments (\S~\ref{subsec:benchmarks} below) suggest that $V_{\max}$ is by far the leading term. Hence, our algorithm yields a significant speedup over the method described in~\cite{AGG08}, as indicated by the ratio $V_{\max}^2/(p \alpha(d))$ of their worst-case time complexities. The same remains true by comparison with the variant of based in the extremality criterion of~\cite{BSS}, 
in which case the factor becomes $2^{d-3}(\log_2V_{\max})^{d-3}/(p\alpha(d))$. 

Observe that, unlike \Call{ComputeExtRays}{}, the other algorithms may temporarily store non-extreme elements of $\CC$ during the inductive step. As a result, the space complexity is not optimal, and it can only be bounded by $O(d \card{V}^2)$. This may be harmful to the scalability of these algorithms.

\begin{table}[t]
\caption{Comparison of the complexity of incremental methods}\label{tab:comparison}
\footnotesize
\begin{center}
\setlength{\extrarowheight}{0.5ex} 
\begin{tabular}{r|@{}c@{}|@{}c@{}|c|@{}c@{}|@{}c@{}|@{}c@{}} 
 \multicolumn{2}{c|}{} & \multirow{2}{*}{$\;$\Call{ComputeExtRays}{}$\;$} & \multirow{2}{2cm}{previous algorithm~\cite{AGG08}} &
\multirow{2}{*}{variant derived from~\cite{BSS}} & 
\multicolumn{2}{c}{classical DDM} \\
\cline{6-7}
\multicolumn{2}{c|}{} & & & &comb. & algebraic \\
\hline 
\multicolumn{2}{c|}{extremality test} & $O(p d \alpha(d))$ & $O(d \card{V}^2)$ &--- & $O(p \card{V})$ & $O(p d^2)$ \\
\hline
\multirow{2}{1cm}{induct. step} & time & $O(p d \alpha(d)\card{V}^2)$ & $O(d \card{V}^4)$ & $O(d\card{V}^22^{d-3}(\log_2\card{V})^{d-3})$ &$O(p \card{V}^3)$ & $O(p d^2 \card{V}^2)$ \\ 
\cline{2-7}
& space & optimal & not optimal & not optimal &\multicolumn{2}{c}{optimal} \\ \hline 
\multicolumn{2}{c|}{speedup (ratio)} & $1$ & $\dfrac{V_{\max}^2}{p \alpha(d)}$ & $2^{d-3}(\log_2V_{\max})^{d-3}/(p\alpha(d))$ & 
$\dfrac{V_{\max}}{d \alpha(d)}$ & $\dfrac{d}{\alpha(d)}$ \\ 
\hline
\multicolumn{2}{c|}{overall} & $O(p^2 d \alpha(d) V^2_{\max})$ & $O(p d V^4_{\max})$ & $\;O(p dV_{\max}^22^{d-3}(\log_2V_{\max})^{d-3})\;$ 
&$\;O(p^2 V^3_{\max})\;$ & $\;O((p d)^2 V^2_{\max})\;$
\end{tabular}
\end{center}
\end{table}

For the sake of completeness, we provide in Table~\ref{tab:comparison} a comparison with the classical double description method, whose principle is close to our algorithm. In the classical case, the elimination of redundant elements can be performed using either an algebraic criterion which can be checked in $O(p d^2)$ arithmetical operations, or a combinatorial criterion of complexity $O(p \card{V})$. See~\cite{FukudaProdon96} for a detailed presentation. We observe that the complexity of the extremality test and the inductive step of our algorithm, as functions of the size of their inputs ($d$, $p$, and $\card{V}$), is smaller than the complexity of their classical analogues (in general, the size of $V$, both in the classical and the tropical settings, is much larger than $d \alpha(d)$).

\subsection{Tropical extreme rays versus tropical polyhedral complexes}

\begin{figure}
\begin{center}
\begin{tikzpicture}[>=triangle 45,shift={(-1,1)},scale=0.7,small arrow/.style={decoration={markings,mark=at position 1 with {\arrow[scale=0.6,#1]{>}}},postaction={decorate},shorten >=0.4pt}]
\node[inner sep=0] at (0,0) {\includegraphics[scale=0.09,viewport=500 400 3500 2400,clip]{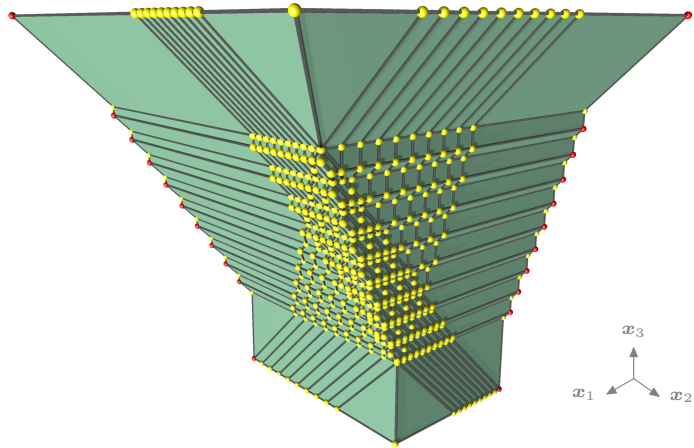}};
\begin{scope}[xshift=5.6cm,yshift=-3.1cm]
\draw[gray,thin,small arrow] (0,0) -- (xyz polar cs:angle=-35,radius=0.6) node[right,font=\tiny] {$\vect{x}_2$};
\draw[gray,small arrow] (0,0) -- (xyz cs:y=0.6) node[above,font=\tiny] {$\vect{x}_3$};
\draw[gray,small arrow] (0,0) -- (xyz polar cs:angle=-150,radius=0.6) node[left,font=\tiny] {$\vect{x}_1$};
\end{scope}
\end{tikzpicture}
\end{center}
\caption{An intersection $10$ affine half-spaces in dimension $3$ which has $24$ vertices (in red) and $1381$ pseudo-vertices (consisting of the $24$ vertices, and the points depicted in yellow)}\label{fig:tropical_poly3}
\end{figure}

Another approach, along the lines of~\cite{DS,Joswig2008}, would consist in representing tropical polyhedra by polyhedral complexes in the usual sense. 
However, an inconvenient of polyhedral complexes is that their number of vertices (called ``pseudo-vertices'' to avoid ambiguities) is exponential in the number of extreme rays~\cite{DS}. Hence, the representations used here are more concise. This is illustrated in Figure~\ref{fig:tropical_poly3} (generated using {\sc Polymake}~\cite{polymake}), which shows an intersection of 10 tropical half-spaces, corresponding to the ``natural'' pattern studied in~\cite{AllamigeonGaubertKatzJCTA2011}. There are only 24 extreme rays, but 1381 pseudo-vertices. 

Another approach, developed by Lorenzo and de la Puente in~\cite{LorenzoPuente2011} (see also Truffet~\cite{Truffet2010}), relies on a similar decomposition of tropical polyhedral cones as polyhedral complexes.
Unlike the previous one, the decomposition arises from half-spaces. Given a cone $\CC$ defined by a system $A \vect{x} \leq B \vect{x}$ with $A, B \in \maxplus^{p \times d}$, the algorithm of~\cite{LorenzoPuente2011} consists in the enumeration of tuples of~$2d$ integers in $\{1, \dots, p\}$ corresponding to (possibly empty) cells of the associated complex. Its worst-case complexity is in $O(p^{2d})$, which is greater than the bound given in~\eqref{eq:complexity_bound} on the complexity of~\Call{ComputeExtRays}{}, as soon as $d \gg 1$ or $p \gg 1$. On top of that, it is possible to exhibit tropical cones with a polynomial number of extreme rays, but on which the algorithm of~\cite{LorenzoPuente2011} have an exponential complexity. Given $p$ reals $t_1 < \ldots < t_p$, consider the polyhedral cone $\CC$ (see~\cite{AllamigeonGaubertKatzJCTA2011}) defined by the inequalities 
$t_k \vect{x}_2 \leq \mpplus_{j \neq 2} t_k^{j-1} \vect{x}_j$, 
for $1 \leq k \leq p$. 
This cone arises as the second polar of the cyclic cone studied in~\cite{AllamigeonGaubertKatzLAA2011}. Here, $t_k^{j-1}$ refers to the scalar $(j-1) \times t_k$. It can be verified that the algorithm of~\cite{LorenzoPuente2011} enumerates $\binom{2p+d-2}{2p}$ tuples. In contrast, as shown in the proof of~\cite[Proposition~6]{AllamigeonGaubertKatzLAA2011}, the number of extreme rays of the cone $\CC$, and also of all intermediate cones $\CC_j$ ($1 \leq j \leq p$) defined by the $j$ first inequalities, is bounded by $O(p d)$. The complexity of the algorithm~\Call{ComputeExtRays}{} is thus polynomial, bounded by $O(p^4 d^3 \alpha(d))$.

\subsection{Benchmarks}\label{subsec:benchmarks}

The algorithm~\Call{ComputeExtRays}{} has been implemented by Allamigeon in the library {\tt TPLib} (Tropical Polyhedral Library) written in OCaml~\cite{TPLib}. Table~\ref{tab:experiments} reports some experiments for different classes of tropical cones: 
\begin{inparaenum}[(i)]
\item samples formed by several cones chosen randomly (referred to as rnd$x$ where $x$ is the size of the sample),
\item and the polars of \emph{signed cyclic cones} which are known to have a very large number of extreme elements~\cite{AllamigeonGaubertKatzJCTA2011}. 
\end{inparaenum} 
For each cone, the first columns respectively report the dimension $d$, the number of half-spaces $p$, the size of the final set of extreme rays, the mean size of the intermediate sets, and the execution time $T$ (for samples of ``random'' cones, we give average results). 

\begin{table}[t]
\caption{Execution time benchmarks of \texttt{TPLib} on a single core of a $3 \, \, \mathrm{GHz}$ Intel Xeon with $3 \, \, \mathrm{Gb}$ RAM}\label{tab:experiments}
\begin{center}
\begin{small}
\begin{tabular}{|c|c|c|c|c|c|c|c|}
\hline
 & $d$ & $p$ & \# final & \# inter. & $T$ (s) & $T'$ (s) & $T/T'$\\
\hline
rnd$100$ & $12$ & $15$ & $32$ & $59$ & $0.24$ & $6.72$ & $0.035$ \\
rnd$100$ & $15$ & $10$ & $555$ & $292$ & $2.87$ & $321.78$ & $8.9 \cdot 10^{-3}$ \\
rnd$100$ & $15$ & $18$ & $152$ & $211$ & $6.26$ & $899.21$ & $7.0 \cdot 10^{-3}$ \\
rnd$30$ & $17$ & $10$ & $1484$ & $627$ & $15.2$ & $4667.9$ & $3.3 \cdot 10^{-3}$ \\
rnd$10$ & $20$ & $8$ & $5153$ & $1273$ & $49.8$ & $50941.9$ & $9.7 \cdot 10^{-4}$ \\
rnd$10$ & $25$ & $5$ & $3999$ & $808$ & $9.9$ & $12177.0$ & $8.1 \cdot 10^{-4}$ \\
rnd$10$ & $25$ & $10$ & $32699$ & $6670$ & $3015.7$ & --- & --- \\
cyclic & $10$ & $20$ & $3296$ & $887$ & $25.8$ & $4957.1$ & $5.2 \cdot 10^{-3}$ \\
cyclic & $15$ & $7$ & $2640$ & $740$ & $8.1$ & $1672.2$ & $5.2 \cdot 10^{-3}$ \\
cyclic & $17$ & $8$ & $4895$ & $1589$ & $44.8$ & $25861.1$ & $1.7 \cdot 10^{-3}$ \\
cyclic & $20$ & $8$ & $28028$ & $5101$ & $690$ & $\sim 45 \text{ days}$ & $1.8 \cdot 10^{-4}$ \\
cyclic & $25$ & $5$ & $25025$ & $1983$ & $62.6$ & $\sim 8 \text{ days}$ & $9.1 \cdot 10^{-5}$ \\
cyclic & $30$ & $5$ & $61880$ & $3804$ & $261$ & --- & --- \\
cyclic & $35$ & $5$ & $155040$ & $7695$ & $1232.6$ & --- & --- \\
\hline
\end{tabular}
\end{small}
\end{center}
\end{table}

In our implementation, inequalities are dynamically ordered during the execution: at each step of the induction, the inequality $\vect{a} \vect{x} \leq \vect{b} \vect{x}$ is chosen so as to minimize the number of combinations $(\vect{a}\vect{v}^j)\vect{v}^i \mpplus (\vect{b} \vect{v}^i) \vect{v}^j$. Note that this strategy does not guarantee that the size of the intermediate sets of extreme elements is smaller. However, it reports better results than without ordering.

We compare our algorithm with an implementation of the algorithm of~\cite{toolbox,AGG08} incorporating the optimizations in~\cite{AllamigeonThesis}, whose execution time $T'$ is given in the seventh column. When the number of extreme rays is of order of $10^4$, the second algorithm needs several days to terminate. For instance, its execution have lasted $45$ days on the signed cyclic cone with the parameters $d = 20$ and $p = 8$, while our algorithm \Call{ComputeExtRays}{} has returned in only $690$ seconds. Therefore, for some extreme cases (for instance $d \geq 30$), the comparison could not be made in practice. 

The ratio $T/T'$ illustrates that our algorithm brings a breakthrough in terms of execution time, which confirms the discussion of Section~\ref{subsec:existing_incremental}. Recall that the main difference between the two algorithms lies in the criterion of elimination of non-extreme vectors, which can be evaluated in time $O(p d \alpha(d))$ (Theorem~\ref{th:complexity}) and $O(d \card{V}^2)$ respectively. As indicated by the experiments, the term $\card{V}^2$ is by far dominating, so that our algorithm benefits from relying on an extremality criterion which does not depend on this factor. This shows that the result of Theorem~\ref{th:main} is interesting not only in theory but also in practice.

\bibliographystyle{amsalpha}
\def\cprime{$'$}
\providecommand{\bysame}{\leavevmode\hbox to3em{\hrulefill}\thinspace}
\providecommand{\MR}{\relax\ifhmode\unskip\space\fi MR }
\providecommand{\MRhref}[2]{%
  \href{http://www.ams.org/mathscinet-getitem?mr=#1}{#2}
}
\providecommand{\href}[2]{#2}

\end{document}